\documentclass{amsart}
\usepackage{hyperref,amssymb,amsrefs}
\newtheorem%
{thm}{Theorem}[section]
\newtheorem%
{proposition}[thm]{Proposition}
\newtheorem%
{lemma}[thm]{Lemma}
\newtheorem%
{lemmadef}[thm]{Lemma-Definition}
\newtheorem%
{corollary}[thm]{Corollary}
\newtheorem%
{conjecture}[thm]{Conjecture}

\theoremstyle{definition}

\theoremstyle{remark}
\newtheorem{remark}[thm]{Remark}

\newcommand{\dontprint}[1]{\relax}

\newcommand{\HH}{\mathit{HH}}

\newcommand{\cO}{{\mathcal{O}}}
\newcommand{\cU}{{\mathcal{U}}}
\newcommand{\cD}{{\mathcal{D}}}

\newcommand{\GL}{\mathit{GL}}
\newcommand{\gl}{\mathit{gl}}
\newcommand{\bp}{{\bar\partial}}

\newcommand{\MC}{{\mathrm{MC}}}
\newcommand{\sgn}{{\mathrm{sgn}}}
\newcommand{\E}{\mathcal{E}}
\newcommand{\barpart}{{\bar\partial}}

\newcommand{\DO}{\mathfrak{D}}

\DeclareMathOperator{\Str}{Str}
\DeclareMathOperator{\Tr}{Tr}
\DeclareMathOperator{\tr}{tr}
\newcommand{\abtop}[2]{\genfrac{}{}{0pt}{}{#1}{#2}}
\DeclareMathOperator{\dist}{dist}
\DeclareMathOperator{\supp}{supp}
\newcommand{\R}{\mathbb{R}}
\newcommand{\C}{\mathbb{C}}

\newcommand{\N}{\mathbb{N}}
\title[RRH formula for traces]{A Riemann--Roch--Hirzebruch formula for traces of differential
operators}
\author{Markus Engeli}
\author{Giovanni Felder}
\address{Department of mathematics, ETH Zurich, 8092 Zurich,
Switzerland}
\begin{document}
\begin{abstract}
Let $D$ be a holomorphic differential operator acting on sections of
a holomorphic vector bundle on an $n$-dimensional compact complex
manifold. We prove a formula, conjectured by Feigin and Shoikhet,
giving the Lefschetz number of $D$ as the integral over the manifold
of a differential form. The class of this differential form is
obtained via formal differential geometry from the canonical
generator of the Hochschild cohomology $\HH^{2n}(\cD_n,\cD_n^*)$ of
the algebra of differential operators on a formal neighbourhood of a
point. If $D$ is the identity, the formula reduces to the
Riemann--Roch--Hirzebruch formula.

\bigskip

\noindent{\sc R\'esum\'e.} Soit $D$ un op\'erateur diff\'erentiel
holomorphe op\'erant sur les sections d'un fibr\'e vectoriel
holomorphe sur une vari\'et\'e complexe de dimension $n$. Nous
d\'emontrons une formule, conjectur\'ee par Feigin et Shoikhet,
donnant le nombre de Lefschetz de $D$ comme int\'egrale d'une forme
diff\'erentielle sur la vari\'et\'e. La classe de cette forme
diff\'erentielle est obtenue, via la g\'eom\'etrie diff\'erentielle
formelle, du g\'en\'erateur canonique de la cohomologie de
Hochschild $\HH^{2n}(\cD_n,\cD_n^*)$ de l'alg\`ebre des op\'erateurs
diff\'erentiels sur un entourage formel d'un point. Si $D$ est
l'identit\'e, la formule se r\'eduit \`a la formule de
Riemann--Roch--Hirzebruch.
\end{abstract}
 \maketitle
 \setcounter{tocdepth}{1}
 \tableofcontents
\section{Introduction}
Let $E\to X$ be a holomorphic vector bundle of rank $r$ on a compact
connected complex manifold $X$ of complex dimension $n$. Let $\cD_E$
be the sheaf of holomorphic differential operators acting on
sections of $E$.

Global differential operators $D\in\cD_E(X)=\Gamma(X,\mathcal D_E)$
act on the sheaf cohomology groups $H^j(X,E)$ of $E$ and thus we
have algebra homomorphisms
\[
H^j:\cD_E(X)\to \mathrm{End}(H^j(X,E)).
\]
Since the cohomology of $E$ is finite dimensional, we can consider
the {\em Lefschetz number} (or supertrace) $L\colon\cD_E(X)\to\C$,
\[
D\mapsto L(D)=\sum_{j=0}^n(-1)^j \mathrm{tr}(H^j(D)).
\]
If $D=\mathrm{Id}$ is the identity then $L(\mathrm{Id})$ is the
holomorphic Euler characteristic of $E$; it is given by the
Riemann--Roch--Hirzebruch theorem as the integral over $X$ of a
characteristic class. Our aim is to generalize this formula to the
case of a general differential operator $D$ by writing the Lefschetz
number as the integral over $X$ of a differential form $\chi_0(D)$
whose value at a point $x\in X$ depends on finitely many derivatives
of the coefficients of $D$ at $x$.

The formula for the differential form $\chi_0$ depends on the choice
of a connection on the holomorphic vector bundles $T^{1,0}X$ and $E$
and is similar to the formula written in \cite{FFS} for the
canonical trace of the quantum algebra of functions in deformation
quantization of symplectic manifolds. Its ingredients are the
Hochschild cocycle of \cite{FFS} and formal differential geometry.
Let $\cD_{n,r}=M_r(\cD_n)$ be the algebra of $r$ by $r$ matrices
with coefficients in the algebra of formal differential operators
$\cD_{n}=\C[[y_1,\dots,y_n]][\partial_{y_1},\dots,\partial_{y_n}]$.
By \cite{FT}, the continuous Hochschild cohomology
$\HH^\bullet(\cD_{n,r},\cD_{n,r}^*)$ is one-dimensional,
concentrated in degree $2n$ and is generated by a $2n$-cocycle
$\tau^r_{2n}\colon\cD_{n,r}^{\otimes (2n+1)}\to\mathbb C$ given in
\cite{FFS} by an explicit integral formula. Formal differential
geometry, see \cite{BR}, gives a realization of $\cD_E(X)$ as the
algebra of horizontal sections for a flat connection $\nabla$ on the
bundle of algebras $\hat \cD_E=J_1E\times_G\cD_{n,r}\to X$ with
fibre $\cD_{n,r}$. Here $J_1E\to X$ denotes the extended frame
bundle, whose fibre at $x\in X$ consists of pairs of bases, one of
$T_x^{1,0}X$ and one of $E_x$; it is a principal bundle for the
group $G=\GL_n(\C)\times\GL_r(\C)$. More generally, let $J_pE$ be
the complex manifold of $p$-jets at $0$ of local bundle isomorphisms
$\C^n\times\C^r\to E$. These manifolds come with holomorphic
$G$-equivariant submersions $J_{p+1}\to J_p$ with contractible
fibres. The flat connection depends on the choice (unique up to
homotopy) of a $G$-equivariant section $\phi\colon J_1E\to J_\infty
E=\varprojlim J_pE$. Such sections can be constructed out of
connections on $J_1E$. Upon local trivialization of $J_1E$ the flat
connection has the form $\nabla(\hat D)=d\hat D+[\omega,\hat D]$ for
some 1-form $\omega$ on $X$ with values in the first order
differential operators in $\cD_{n,r}$ and the isomorphism
$\cD_E(X)\to\mathrm\mathrm{Ker}(\nabla)$ sends $D$ to its Taylor
expansion $\hat D=\phi_*D$ with respect to the local coordinates and
trivialization of $E$ given by $\phi$.

With these notations the formula for $\chi_0(D)$ in terms of the
horizontal section $\hat D$ associated with $D$ is
\[
\chi_0(D)=\tau^r_{2n}(\hat D,\omega,\dots,\omega).
\]
The multilinear form $\tau^r_{2n}$ on $\cD_{n,r}$ is extended to
differential forms with values in $\cD_{n,r}$ by linearity: if
$\omega=\sum\omega_jdx_j$ in terms of local real coordinates $x_j$,
$j=1,\dots,2n$,
\[
\chi_0(D)=\sum \tau^r_{2n}(\hat
D,\omega_{j_1},\dots,\omega_{j_{2n}})\,dx_{j_1}\wedge\cdots\wedge
dx_{j_{2n}}.
\]
The local objects $\hat D$ and $\omega$ depend on a choice of a
local trivialization of $J_1E$, but the differential form $\chi_0$
is globally defined as a consequence of the fact that $\tau^r_{2n}$
is basic for the action of $G$.  Our main result is
\begin{thm}\label{t-1} For any $D\in\cD_E(X)$,
\[
L(D)=\frac{1}{(2\pi i)^n}\int_X\chi_{0}(D).
\]
\end{thm}
Moreover, for the identity differential operator, it is known
\cites{FT,NT1} that the class of $\chi_0(\mathrm{Id})$ is the
component of degree $2n$ of the Hirzebruch class
$\mathrm{td}(T_X)\mathrm{ch}(E)$ and thus we recover the
Riemann--Roch--Hirzebruch theorem. Also, the direct calculation of
\cite{FFS} shows that $\chi_0(\mathrm{Id})$ is the representative of
the Hirzebruch class given by the Chern--Weil map in terms of the
curvature of the connection on $T^{1,0}X\oplus E$ canonically
associated with $\phi$.

The proof of the theorem is obtained by showing that the linear
functions $T_1=L$ and $T_2=\int_X\chi_0$ on the Hochschild 0-th
homology
\[
\HH_0(\cD_E(X))=\cD_E(X)/[\cD_E(X),\cD_E(X)]
\]
are proportional to a third linear function $T_3$ constructed
essentially in \cites{BG,Wodzicki}: a global differential operator
$D\in\cD_E(X)$ defines a global 0-cycle in the complex of sheaves
$\mathcal C_\bullet(\cD_E)$ of Hochschild chains of $\cD_E$, which
is quasi-isomorphic to the complex of sheaves ${\mathbb C}_X[2n]$ of
locally constant continuous functions concentrated in degree $-2n$.
Thus there is a map $T_3\colon \HH_0(\cD_E(X))\to H^0(X,{\mathbb
C}_X[2n]) = H^{2n}(X,\C)\simeq \C$.

The statement of Theorem \ref{t-1} was conjectured around 2001 by B.
Feigin and B. Shoikhet. In the case of curves a formula for $L(D)$
in terms of residues had been found by A. Beilinson and V.
Schechtman (Lemma 2.2.3 in \cite{BS}, see also \cite{S}). A formula
for the normalized trace in deformation quantization of a symplectic
manifold, analogous to the one of Theorem \ref{t-1} was proposed in
\cite{FFS}. The proof of that formula is simpler since the space of
traces is one-dimensional in that situation, so one just has to
check the normalization. The difficulty here is that
$\HH_0(\cD_E(X))$ is not one-dimensional in general. An indirect
approach to proving that $T_1=T_3$, proposed in \cite{FLS}, is to
embed $\cD_E(X)$ in a suitable complex of algebras with
one-dimensional cohomology and show that both $T_1$ and $T_3$ extend
to chain maps on this complex. If the Euler characteristic of $E$
does not vanish one can then deduce from the classical
Riemann--Roch--Hirzebruch theorem that $T_1=C\cdot T_3$ for some
$C$. The rigorous completion of this programme presents some
technical difficulties but it should lead to a proof of $T_1=T_3$ if
$E$ has non-vanishing Euler characteristic. In a very recent
preprint \cite{R}, A. Ramadoss shows that the approach of \cite{FLS}
could be extended to the much more general case where $X$ admits a
vector bundle with non-vanishing Euler characteristic.

Our result gives in particular a different direct proof of the fact
that $T_1=T_3$, without assumptions on $X$ or $E$. It does not use
the Riemann--Roch--Hirzebruch theorem. 

\medskip

\noindent{\bf Acknowledgements.} This work has been partially
supported by the European Union through the
FP6 Marie Curie RTN ENIGMA (Contract number MRTN-CT-2004-5652), and the
Swiss National Science Foundation (grants 21-65213.01 and 200020-105450).

We are grateful to Boris Shoikhet for introducing us to this subject
and for sharing his insights, which were essential to this project.
We also thank Vasilyi Dolgushev, Pavel Etingof, Victor Ginzburg,
David Kazhdan and Carlo Rossi for discussions and suggestions. We
thank the referee for carefully reading the manuscript and

suggesting several improvements.

\section{Hochschild homology of the algebra of differential operators}
\subsection{Hochschild homology}\label{ss-Piscine Molitor}
Let $A$ be an algebra over $\mathbb C$ with unit 1 and set $\bar
A=A/\mathbb C1$. We denote $\bar a$ the class in $\bar A$ of $a\in
A$. The Hochschild homology $\HH_\bullet(A)$ of $A$ with
coefficients in the bimodule $A$ is the homology of the (normalized)
Hochschild chain complex $\cdots\stackrel b\to C_q(A)\stackrel b\to
C_{q-1}(A)\stackrel b\to\cdots$ with
\[
C_q(A)=A\otimes \bar A^{\otimes q}, \qquad q\geq 0,
\]
and differential
\begin{eqnarray}\label{e-Hochschild}
b(a_0,\dots,a_q)&=&\sum_{j=0}^{q-1}(-1)^j(a_0,\dots,a_ja_{j+1},\dots,a_q)
\\
&&+(-1)^{q}(a_qa_0,a_1,\dots,a_{q-1}).\notag
\end{eqnarray}
Here $a_0,\dots,a_q\in A$ and we write $(a_0,\dots,a_q)$ instead of
$a_0\otimes \bar{a}_1\otimes\cdots\otimes\bar{a}_q$. For topological
algebras one has to take the projective tensor product, as explained
in \cite{Connes}, Ch.~II.

 Let $\cO_n=\C[[y_1,\dots, y_n]]$ be the algebra of
formal powers series in $n$ variables and
$\cD_n=\cO_n[\partial_{y_1},\dots,\partial_{y_n}]$ the algebra of
formal differential operators. Let also
$\cO_n^{\mathrm{pol}}=\mathbb C[y_1,\dots,y_n]$,
$\cD_n^{\mathrm{pol}}=\cO_n^{\mathrm{pol}}[\partial_{y_1},\dots,\partial_{y_n}]$
be the subalgebras of polynomial functions and differential
operators. As shown by Feigin and Tsygan \cite{FT}, the Hochschild
homology of $\cD_n^\mathrm{pol}$ is one-dimensional and concentrated
in degree $2n$. A representative of a generator of
$\HH_{2n}(\cD_{n})$ in the normalized Hochschild chain complex is
\[
c_{2n}=\sum_{\pi\in S_{2n}}\sgn(\pi)\,1\otimes
u_{\pi(1)}\otimes\cdots\otimes u_{\pi(2n)},\quad
u_{2j-1}=\partial_{y_j},\,u_{2j}=y_j.
\]
Thus there is a unique linear form on Hochschild homology whose
value on $c_{2n}$ is one. This linear form is the class of a cocycle
in the complex dual to the Hochschild complex. An explicit formula
for such a cocycle $\tau_{2n}$ was found in \cite{FFS}. It has the
following properties.
\begin{enumerate}
\item[(i)] $\tau_{2n}$ extends to a linear form on
$\cD_n^{\otimes(2n+1)}$ obeying the cocycle condition
$\tau_{2n}\circ b=0$, where $b$ is the Hochschild differential, and
the normalization condition: $\tau_{2n}(D_0,\dots,D_{2n})=0$ if
$D_j=1$ for some $j\geq1$.
\item[(ii)] $\tau_{2n}$ is invariant under the action of $\GL_n(\C)$
on $\cD_n$ by linear coordinate transformations. Moreover, if $a=\sum
a_{jk}y_k\partial_{y_j}+b$,
$a_{jk},b\in\mathbb C$, then \\
$\sum_{j=1}^{2n}(-1)^j\tau_{2n}(D_0,\dots,D_{j-1},a,D_{j},\dots,D_{2n-1})=0$.
\item[(iii)] $\tau_{2n}(c_{2n})=1$.
\end{enumerate}
More generally, let $M_r(A)\simeq M_r(\mathbb C)\otimes A$ denote
the algebra of $r$ by $r$ matrices with entries in an associative
algebra $A$.
 Since Hochschild homology is Morita invariant,
$\HH_\bullet(M_r(\cD_n))\simeq \HH_\bullet(\cD_n)$ is also one-dimensional
and is spanned by $c_{2n}$ where we view $\cD_n$ as a subalgebra of
$M_r(\cD_n)$ via $D\to \mathrm{Id}\otimes D$. Define a cocycle
$\tau^r_{2n}$ by
\[
\tau^r_{2n}(A_0\otimes D_0,\dots,A_{2n}\otimes D_{2n})=
\mathrm{tr}(A_0\cdots A_{2n})\tau_{2n}(D_0,\dots,D_{2n}),
\]
$A_i\in M_r(\mathbb C)$, $D_i\in \cD_n$. As a consequence of the
properties of $\tau_{2n}$, $\tau^r_{2n}$ obeys:
\begin{enumerate}
\item[(i)] $\tau^r_{2n}$ is a linear form on
$M_r(\cD_n)^{\otimes(2n+1)}$ obeying the cocycle condition
$\tau^r_{2n}\circ b=0$ and $\tau^r_{2n}(D_0,\dots,D_{2n})=0$ if, for
some $j\geq1$, $D_j$ is the multiplication by a constant matrix.
\item[(ii)] $\tau^r_{2n}$ is invariant under the action of $G=\GL_n(\C)\times\GL_r(\C)$
where $\GL_r(\C)$ acts on $M_r(\cD_{n})$ by conjugation. Moreover, if $a=\sum
a_{jk}y_k\partial_{y_j}+b$, $a_{jk}\in\mathbb C$, $b\in M_r(\mathbb
C)$ then
\[
\sum_{j=1}^{2n}(-1)^j\tau^r_{2n}(D_0,\dots,D_{j-1},a,D_{j},\dots,D_{2n-1})=0.
\]
\item[(iii)] $\tau^r_{2n}(c_{2n})=r$.
\end{enumerate}

\begin{remark}
For any associative algebra $A$ denote by $A_{\mathrm{Lie}}$ the Lie
algebra $A$ with bracket $[a,b]=ab-ba$. Then $A_{\mathrm{Lie}}$ acts
on $C_p(A)$ via
\[
L_a(a_0,\dots,a_p)=\sum_{j=0}^p(a_0,\dots,[a,a_j],\dots,a_p),\qquad
a\in A_{\mathrm{Lie}}
\]
and we have a Cartan formula $L_a=b\circ
\iota_a+\iota_a\circ b$ with
\[
\iota_a(a_0,\dots,a_p)=\sum_{j=1}^p(-1)^{j+1}(a_0,\dots,a_{j-1},a,a_j,\dots,a_{p}).
\]
It follows that $A_{\mathrm{Lie}}$ acts trivially on the cohomology.
The property (ii) may be rephrased as saying that $\tau_{2n}$ is
$G$-{\em basic}, namely $G$-invariant and obeying
$\tau_{2n}^r\circ\iota_a=0$, for $a$ in the Lie algebra of $G$
embedded in $\cD_{n,r}$ as a Lie algebra of first order operators.

It also follows that the cohomology class of $\tau^r_{2n}$ is
invariant under coordinate transformations.
\end{remark}

\subsection{Hochschild chain complex of the sheaf of differential operators}
Let $\cD_E$ be the sheaf of differential operators on $E$. In terms
of holomorphic coordinates and a local holomorphic trivialization of
$E$, a local section of $\cD_E$ has the form
\[
\sum_{I}
a_I(z_1,\dots,z_n)\partial_{z_1}^{i_1}\cdots\partial_{z_n}^{i_n},\qquad
I=(i_1,\dots,i_n)\in\mathbb Z_{\geq0}^n,
\]
with holomorphic matrix-valued coefficients $a_I$, vanishing except
for finitely many multiindices $I$. The sheaf $\cD_E$ is a sheaf of
locally convex algebras: for any open set $U\subset X$, the locally
convex subalgebra $\cD_E(U)^{\leq k}$ of operators of order at most
$k$ is the space of sections of some vector bundle over $U$ and has
the topology of uniform convergence on compact subsets. Then the
inductive limit $\cD_E(U)=\cup_{k}\cD_E(U)^{\leq k}$ with the
inductive limit topology is a complete locally convex algebra. This
is the topology considered in \cite{BG}. Then one has the following
result:

\begin{thm}\label{l-Miss Kenton} \cites{BG,Wodzicki} Every point of $X$ has
a coordinate neighbourhood $U$ such that $\HH_{p}(\cD_E(U))=0$ for
$p\neq 2n$ and $\HH_{2n}(\cD_E(U))$ is one-dimensional generated by
the class of
\[
c_E(U)=\sum_{\pi\in
S_{2n}}\sgn(\pi)(1,x_{\pi(1)},\dots,x_{\pi(2n)}),
\]
where $x_{2j-1}=\partial_{z_j},x_{2j}=z_j$. Here we identify $x\in
\cD(U)$ with the multiple of the identity $\mathrm{Id}_r\otimes x\in
M_r\otimes\cD(U)\simeq \cD_E(U)$, with respect to some
trivialization of $E$.
\end{thm}

\subsection{Formal differential geometry}
We recall some notions of formal differential geometry
\cites{G,GK,GKF}, following \cite{BR}.

Let $W_n=\oplus_i\cO_n\partial_{y_i}$ be the Lie algebra of formal
vector fields and $\gl_r(\cO_n)$ denote $M_r(\cO_n)$ viewed as a Lie
algebra, with commutator bracket. The Lie algebra $W_n$ acts on
$\gl_r(\cO_n)$ by derivations and we can thus define the semidirect
product
\[
W_{n,r}=W_n\ltimes\gl_r(\cO_n).
\]
This Lie algebra is embedded in $M_{r}(\cD_n)$ (viewed as Lie
algebra with commutator bracket) as a Lie subalgebra of first order
differential operators. It should be regarded as the Lie algebra of
infinitesimal automorphisms of the trivial bundle of rank $r$ over a
formal neighbourhood of $0\in\C^n$.

A {\em local parametrization} of $E$ is a holomorphic bundle
isomorphism $U\times \C^r\to E|_V$ from the trivial bundle over some
neighbourhood $U\subset\C^n$ of $0$ to the restriction of $E$ to
some open set $V$. Let $J_p E$ be the complex manifold of $p$-jets
at $0\in\C^n$ of local parametrizations. In particular, $J_1 E$ is
the extended frame bundle, whose fibre at $x\in X$ is the space of
pairs of bases of the holomorphic tangent space at $x$ and the fibre
of $E$ at $x$ respectively. The group $G=\GL_n(\C)\times \GL_r(\C)$
acts freely on the right on each $J_p E$, $p=1,2,\dots$ by linear
transformations of $\C^n\times \C^r$ and $J_1E$ is a principal
$G$-bundle over $X$. The complex manifolds $J_pE$ form a projective
system with surjective $G$-equivariant submersions $J_p E\to J_q E$,
$p>q$. The projective limit $J_\infty E$ is, in the language of
\cite{BR}, a holomorphic {\em principal $W_{n,r}$-space}. Namely,
there is a Lie algebra homomorphism $W_{n,r}\to \mathcal V(J_\infty
E)$ from $W_{n,r}$ to the Lie algebra of holomorphic vector fields
on $J_\infty E$, which is an isomorphism $W_{n,r}\to T^{1,0}_\phi
J_\infty E$ at each point $\phi\in J_\infty E$. The inverse map
defines a holomorphic one-form $\Omega_{\MC}\in\Omega^{1,0}(J_\infty
E,W_{n,r})$ with values in $W_{n,r}$ and the homomorphism property
is equivalent to the Maurer--Cartan equation
\[
d\Omega_{\MC}+\frac12[\Omega_{\MC},\Omega_{\MC}]=0.
\]
Moreover, the fibres of the bundle $J_\infty E/G\to J_1E/G=X$ are
contractible and therefore there exists a smooth section (unique up
to homotopy) $\phi:X\to J_\infty E/G$ or, equivalently, a smooth
$G$-equivariant section $\tilde \phi:J_1E\to J_\infty E$. The
Maurer--Cartan form $\Omega_{\MC}$ pulls back to a $G$-equivariant
$1$-form $\tilde\phi^*\Omega_{\MC}$ on $J_1E$ obeying the
Maurer--Cartan equation. This induces a flat connection on the
associated bundle
\[
\hat\cD_E=J_1E\times_G M_r(\cD_n)\to X.
\]
The horizontal sections are in one-to-one correspondence with global
differential operators: to $D\in\cD_E(X)$ there corresponds the
horizontal section $\hat D$. Its value at $x\in X$ is the Taylor
expansion at $0$ of $D$ with respect to the coordinates and the
trivialization defined by $\phi$ at the point $x$. Conversely, every
horizontal section comes from a differential operator. In explicit
terms, let us choose a local trivialization of $J_1E=U\times G$ over
$U\subset X$. Then the  restriction of $\phi$ to $U$ is given by a
map $\phi^U\colon U\to J_\infty E|_U$ and
$\omega=\phi^*\Omega_{\MC}$ is a $W_{n,r}$-valued 1-form on $U$. The
Taylor expansion $\hat D$ is given on $U$ by a map $U\to
M_r(\cD_n)$, $x\mapsto \hat D_x$ obeying
\[
d \hat D+[\omega, \hat D]=0.
\]
A change of trivialization is given by a gauge transformation
$g\colon U\to G$. The section changes as $\hat D_x\mapsto g_x\cdot
\hat D_x$ and $\omega$ as $\omega\mapsto g\cdot\omega-dg g^{-1}$ and
$dg g^{-1}$ is a 1-form with values in the Lie algebra of $G$,
embedded in $M_r(\cD_n)$ as the Lie algebra of first order operators
of the form $\sum a_{jk}y_k\frac\partial{\partial y_j}+b$,
$a_{jk}\in\mathbb C$, $b\in M_r(\mathbb C)$.

\begin{proposition}\label{p-Naoko}
Let $\Omega^\bullet$ be the complex of sheaves of complex-valued
smooth differential forms on $X$ with de Rham differential and let
$\mathcal C(\cD_E)$ be the complex of sheaves of Hochschild chains
of $\cD_E$. There is a homomorphism of complexes of sheaves
\[
\chi_\bullet\colon\mathcal{C}_\bullet(\cD_E)\to \Omega^{2n-\bullet},
\]
depending on a choice of section of $J_\infty E/G\to X$, inducing an
isomorphism of the cohomology sheaves. The map
$\chi_0\colon\cD_E(X)\to\Omega^{2n}(X)$ on global differential
operators is the map appearing in Theorem \ref{t-1} and $\chi_{2n}$
 maps $(D_0,\dots,D_{2n})\in \mathcal{C}_{2n}(U)$ to the function
$\tau^r_{2n}(\hat D_0,\dots,\hat D_{2n})$ on the open set $U$.
\end{proposition}

The rest of this section is dedicated to the construction of
$\chi_\bullet$ and the proof of Proposition \ref{p-Naoko}.

\subsection{Shift by a Maurer--Cartan element}
We start by a general construction on the chain complex of an
arbitrary differential graded algebra. Let $A=\oplus_{j\in\mathbb Z}
A^j$ be a differential graded algebra with unit and with
differential $d\colon A^{j}\to A^{j+1}$. We denote by $|a|=j$ the
degree of a homogeneous element $a\in A^j$. The Hochschild chain
complex of $A$ is $C^\bullet(A)=\oplus_{p\in\mathbb Z}C^p(A)$ with
\[C^p(A)=\Pi_{\sum j_r-q=p}A^{j_0}\otimes \bar A^{j_1}\otimes\cdots\otimes
\bar A^{j_q}.\]
 The differential of the Hochschild complex is
defined as the total differential $\delta=b+(-1)^pd$ on 
$C^p(A)$.\footnote{We introduce the upper index notation $C^q$ to
have a differential of degree one as $d$ is. Thus in the ungraded
case we have $C^q(A)=C_{-q}(A)$, concentrated in negative degrees.}
The Hochschild differential $b$ is defined as in
\eqref{e-Hochschild} except that the last term has an additional
sign $(-1)^{|a_p|(|a_0|+\cdots+|a_{p-1}|)}$ and the differential $d$
is extended as a derivation of degree 1 for the tensor product:
\[
d(a_0,\dots,a_p)=\sum_{j=0}^p(-1)^{|a_0|+\cdots+|a_{j-1}|}(a_0,\dots,da_j,\dots,a_p).
\]
A {\em Maurer--Cartan element} of $A$ is an element $\omega\in A^1$
of degree 1 obeying the Maurer--Cartan equation
\[
d\omega+\omega^2=0
\]
The Maurer--Cartan equation implies that the linear endomorphism
$d_\omega$ of $A$ given by $d_\omega a=da+\omega
a-(-1)^{|a|}a\omega$ is a differential. Moreover $d_\omega$ is a
derivation of degree 1 of the algebra $A$ and therefore the algebra
$A$ with differential $d_\omega$ is a differential graded algebra.
We call this differential graded algebra the twist of $A$ by
$\omega$ and denote it $A_\omega$.

The symmetric group $S_p$ acts on $A\otimes \bar A^{p}$ by
permutations of the last $p$ factors with signs: the transposition
of neighbouring factors $a$ and $b$ is accompanied by the sign
$(-1)^{|a|\cdot|b|}$. Recall that the shuffle product $C_p(A)\otimes
C_q(A)\to C_{p+q}(A)$ is defined by
\[
(a_0,\dots,a_p)\times(b_0,\dots, b_q)=(-1)^{|b_0|\sum
|a_j|}\mathrm{sh}_{p,q} (a_0b_0,a_1,\dots,a_p,b_1,\dots,b_q),
\]
where $\mathrm{sh}_{p,q}x=\sum_{\pi\in S_{p,q}}\sgn(\pi)\pi\cdot x$,
with sum over $(p,q)$-shuffles in $S_{p+q}$, namely over the
permutations that preserve the ordering of the first $p$ and of the
last $q$ letters. The shuffle product is associative and if $A$ is
Abelian (which we do not assume) it is a homomorphism of complexes,
see \cites{Maclane,Loday}.

\begin{proposition}\label{p-Grace}
Let $A_\omega$ be the twist of $A$ by a Maurer--Cartan element
$\omega\in A^1$. Let $(\omega)_k=(1,\omega,\dots,\omega)$ with $k$
factors of $\omega$. Then the map
\[
(a_0,\dots,a_p)\to\sum_{k\geq
0}(-1)^k(a_0,\dots,a_p)\times(\omega)_k
\]
is an isomorphism of complexes $C(A_\omega)\to C(A)$.
\end{proposition}

We split the proof into a few steps.
\begin{lemma} $b(\omega)_0=0$ and, for $k\geq 1$,
$ b(\omega)_k=d(\omega)_{k-1}. $
\end{lemma}
\begin{proof}
The first statement is obvious. Let $k\geq 1$. Then
\begin{eqnarray*}
b(\omega)_k&=&b(1,\omega,\dots,\omega)\\
&=&(\omega,\dots,\omega)+\sum_{j=1}^{k-1}(-1)^j(1,\omega,\dots,{\omega^2},\dots,\omega)
\\
&&+(-1)^k(-1)^{k-1}(\omega,\dots,\omega)\\
\\
&=& -\sum_{j=1}^{k-1}(-1)^j(1,\omega,\dots,{d\omega},\dots,\omega)
\\
&=&d(\omega)_{k-1}.
\end{eqnarray*}
\end{proof}

\begin{lemma} Let 
 $a\in C^p(A)$.Then
\begin{eqnarray*} 
b\left(a\times(\omega)_k\right)&=&b\,
a\times(\omega)_k+(-1)^pa\times
b(\omega)_k\\
&&-(-1)^p\sum_{k=0}^p(-1)^{|a_0|+\cdots+|a_{k-1}|}(a_0,\dots,[\omega,a_k],\dots,a_p)\times(\omega)_{k-1},
\end{eqnarray*}
where $[a,a']=aa'-(-1)^{|a|\cdot|a'|}a'a$ is the graded commutator.
\end{lemma}

\begin{proof} For simplicity, we give the proof in the case where all $a_j$ are of
degree $0$, which is the case appearing in our application. The
additional signs appearing in the general case can be treated
easily.

If we write out the sum over shuffles we see that there are four
types of terms appearing on the left-hand side: those containing the
products $a_ja_{j+1}$, $\omega^2$, $\omega a_j$ and $a_j\omega$. The
terms of the first and of the second type combine to give the first
two terms on the right-hand side. The proof that the signs match is
the same as in the proof of the homomorphism property for
commutative algebras, see \cite{Loday}, Proposition 4.2.2, so we
consider only the last two types. Consider a shuffle $\pi$ appearing
on the left-hand side such that ${l}$ out of the $k$ factors
$\omega$ have been shuffled to the left of $a_j$. Then the term
containing the product $\omega a_j$ comes with a sign
$\sgn(\pi)(-1)^{j-1+{l}}$. The same term occurs for a shuffle $\pi'$
in $(a_0,\dots,[\omega,a_j],\dots,a_p)\times(\omega)_{k-1}$ with a
sign equal to $\sgn(\pi')(-1)^{{l}-1}$, where $(-1)^{{l}-1}$ is the
Koszul sign coming from the fact that ${l}-1$ factors $\omega$ are
permuted by $\pi'$ to the left of the odd element $[\omega,a_j]$.
The signs of the shuffles are related by
$\sgn(\pi)=\sgn(\pi')(-1)^{p-j+1}$. The ratio of signs is thus
$(-1)^{p-1}$, as claimed. The same reasoning can be applied to
$a_j\omega$.
\end{proof}

\begin{lemma}\label{Ian} 
 Let $a\in C^p(A)$ and set $\delta_\omega
a=b\,a+(-1)^pd_\omega a$. Then
\[
\delta\sum_{k\geq 0}(-1)^k a\times(\omega)_k=
\sum_{k\geq0}(-1)^k\delta_\omega a\times(\omega)_k
\]
\end{lemma}
\begin{proof}
This follows from the previous Lemma by inserting the definitions
and summing over $k$. \end{proof}
\begin{lemma}\label{Isabel}
The map $C(A_\omega)\to C(A)$ of Proposition \ref{p-Grace} is an
isomorphism.
\end{lemma}
\begin{proof}
The map is the shuffle multiplication by $\psi=\sum
(-1)^k(\omega)_k$. We claim that the inverse map is the shuffle
multiplication by $\bar\psi=\sum(\omega)_k$. To prove this, we use
that the shuffle product is associative, so that it suffices to show
that $\psi\times\bar\psi=1$. This follows from
\[
(\omega)_k\times(\omega)_{l}=\sum_{\pi\in
S_{k,{l}}}(\omega)_{k+{l}}=\left(\begin{array}{c}k+{l}\\
k\end{array}\right)(\omega)_{k+{l}}.
\]
\end{proof}
Proposition \ref{p-Grace} follows from the last two Lemmata.
\subsection{Hochschild and de Rham cohomology}
We construct a homomorphism of complexes of sheaves
\[ \chi_\bullet\colon\mathcal C_\bullet(\cD_E)\to\Omega^{2n-\bullet}\]
from the Hochschild chain complex of the sheaf $\cD_E$ to the sheaf
of smooth de Rham forms. It is based on formal geometry and thus
depends on a choice of section of $J_\infty E/G$ (but the map
induced on homology is canonical). The map $\chi_0\colon
\cD_E(X)\to\Omega^{2n}(X)$ on global differential operators is the
one appearing in Theorem \ref{t-1}.

To do this we apply the previous constructions to the smooth de Rham
complex $A=\Omega(U,\hat\cD_E)$ with values in the vector bundle
$\hat\cD_E$ on some open subset $U\subset X$. Let $\hat D\in A^0$
denote the horizontal section corresponding to a differential
operator $D\in\cD_E(U)$. Locally, upon trivialization of $T^{1,0}X$
and $E$, the condition of horizontality is $d \hat D+[\omega,\hat
D]=0$ for some Maurer--Cartan element $\omega$.
\begin{proposition}\label{p-Asfar}
Let $U$ be a sufficiently small open neighbourhood of
any point in $X$. Let $D_0,\dots,D_p\in \cD_E(U)$ be differential
operators on $U$ and $\hat D_0,\dots,\hat D_p\in A^0$ be the
corresponding horizontal sections of $\hat\cD_E$ on $U$. Then the differential
$(2n-p)$-forms on $U$ 
\begin{equation}\label{e-Stevens}
 \chi_p(D_0,\dots,D_p)=
 \tau^r_{2n}\left(\mathrm{sh}_{p,2n-p}(\hat
 D_0,\hat D_1,\dots,\hat D_p,\omega,\dots,\omega)\right),
\end{equation}
are well-defined (i.e., independent of the trivialization of
$J_1E$), continuous, and obey the relations
\[ 
d\circ\chi_{p}=(-1)^{p-1}\chi_{p-1}\circ b.
\]
\end{proposition}

\begin{proof}
If we change trivialization of the extended frame bundle $J_1E$,
then $\hat D$, $\omega$ change by the action of an element of $G$,
under which $\tau^r_{2n}$ is invariant, and the shift of $\omega$ by
a one-form with values in the Lie algebra of $G$ embedded in $A$. By
property (ii) of $\tau_{2n}^r$, see \ref{ss-Piscine Molitor}, the
right-hand side of \eqref{e-Stevens} is unaffected by such a shift.
The continuity is clear: since $\tau_{2n}^r$ depends non-trivially
only on finitely many Taylor coefficients of its arguments, the
$C^\ell$-norms on compact subsets of $\chi_p(D_0,\dots,D_p)$ are
estimated by $C^{\ell'}$-norms of the coefficients of $D_0,\dots,
D_p$ which by analyticity are in turn controlled by sup norms on
(slightly larger) compact subsets.

In the notation of Proposition \ref{p-Grace},
\[ 
 \chi_p(D_0,\dots,D_p)=
 \tau^r_{2n}\left((\hat D_0,\dots,\hat
 D_p)\times(\omega)_{2n-p}\right).
\]
The cochain $a=(\hat D_0,\dots,\hat D_p)$ obeys $d_\omega(a)=0$,
thus the homomorphism property of Proposition \ref{p-Grace} reduces
to
\[
\delta \sum_{k\geq 0}(-1)^k a\times(\omega)_k= \sum_{k\geq 0}(-1)^k
b\,a\times(\omega)_k,\qquad \delta=b\pm d.
\]
The component of Hochschild degree $2n$ is
\[ 
(-1)^{p-1}b(a\times(\omega)_{2n-p+1})+(-1)^p(-1)^pd(a\times(\omega)_{2n-p})=
(-1)^{p-1}b\,a\times(\omega)_{2n-p+1}.
\]
If we apply $\tau^r_{2n}$ the first term vanishes and we obtain the
claim.
\end{proof}

Since the expression on the right-hand side of \eqref{e-Stevens} are
local, it is clear that $\chi_p$ are compatible with the inclusion
of open sets and thus define maps of sheaves. Moreover, by the
normalization of $\tau^r_{2n}$, we see that $\chi_{2n}(c_E(U))=r$,
where $c_E(U)$ is the generator of Theorem \ref{l-Miss Kenton}. Thus
$\chi_\bullet$ induces a non-trivial map on homology. By Theorem
\ref{l-Miss Kenton} this map is an isomorphism. This concludes the
proof of Proposition \ref{p-Naoko}.

\section{The third trace}\label{s-Midori}
The idea of the proof of Theorem \ref{t-1} is to show that the two
traces $T_1(D)=L(D)=\sum (-1)^j\mathrm{tr} (H^j(D))$ and
$T_2(D)=\int_X\chi_0(D)$ are both proportional to a third linear
form $T_3\colon \cD_E(X)\to\C$  constructed via Theorem \ref{l-Miss
Kenton} with the help of a finite open cover $\cU=(U_\alpha)$.

Consider the Hochschild complex of sheaves $\mathcal
C_\bullet(\cD_E)$:
\[
\dots\to
\cD_E\otimes\bar\cD_E\otimes\bar\cD_E\to\cD_E\otimes\bar\cD_E\to\cD_E\to
0
\]
Let $\mathcal U=(U_\alpha)$ be a sufficiently fine open cover of
$X$. Let $C^{p,q}=\check C^q(\mathcal U,\mathcal C_{-p}(\cD_E))$,
$(q\geq 0,p\leq0)$, where $\check C^q(\mathcal U,\mathcal
F)=\oplus_{\alpha_0<\dots<\alpha_{p}}\mathcal
F(U_{\alpha_0}\cap\cdots\cap U_{\alpha_p})$ be the Hochschild--\v
Cech double complex. Global differential operators $D\in\cD_E(X)$
define cocycles in $C^{0,0}$. The restriction $D|_{U_\alpha}$ of $D$
to a sufficiently small open set is a Hochschild boundary by Theorem
\ref{l-Miss Kenton}. Thus $D|_{U_\alpha}=b D^{(1)}_\alpha$ for some
$D^{(1)}\in C^{-1,0}$.
\begin{figure}
\begin{tabular}{lllll|l}
  $s$                  &       &                       &    &     & $2n$     \\
 $\uparrow$            &       &                       &    &     &          \\
\!$D^{(2n)}\rightarrow$&       &                       &    &     & $2n-1$   \\
                       &$\ddots$&                      &    &     &          \\
                       &       & $\uparrow$            &    &     &          \\
                       &       & \!$D^{(2)}\rightarrow$&    &     & $1$      \\
                       &       &     &          $\uparrow$  &     &          \\
                       &       &     &\!$D^{(1)}\rightarrow$& $D$ & $0$      \\
\hline
           $-2n$       &       &     &               $-1$   & $0$ &          \\
           \end{tabular}
\caption{The Hochschild--\v Cech double complex}\label{f-one}
\end{figure}
Since $b$ and the \v Cech differential commute, $(\check\delta
D^{(1)})_{\alpha\beta}=D^{(1)}_\beta-D^{(1)}_\alpha$ is a Hochschild
cycle for the algebra $\cD_E(U_\alpha\cap U_\beta)$ and is thus
exact. Continuing in this way we can ``climb the staircase'', see
Fig.~\ref{f-one}, and find $D^{(j)}\in C^{-j,j-1}$, $j=1,\dots,2n$,
such that
\begin{equation}\label{e-staircase}
b D^{(1)}=D,\qquad\check\delta D^{(j)}=b D^{(j+1)},\qquad
j=1,\dots,2n-1.
\end{equation}
Finally we get to the point where the Hochschild homology is
nontrivial and we obtain
\begin{equation}\label{e-Lord Darlington}
\check \delta D^{(2n)}=s+b D^{(2n+1)},
\end{equation}
where $s\in C^{2n,-2n}$ has the form
\begin{equation}\label{e-Toru Watanabe}
s_{\alpha_0,\dots,\alpha_{2n}}=\lambda_{\alpha_0,\dots,\alpha_{2n}}(D)
c_{E}(U_{\alpha_0}\cap\cdots\cap U_{\alpha_{2n}}),
\end{equation}
for some \v Cech cocycle  $\lambda(D)\in \check C^{2n}(\mathcal
U,\mathbb{C})$ with values in the locally constant sheaf $\mathbb
C_X$. Its class $[\lambda(D)]\in H^{2n}(X,\mathbb C)\simeq \C$ is
(up to sign) $T_3(D)$.

\section{The first trace is proportional to the third\dots}\label{s-1=3}
Here we study the first trace $T_1=L$ and describe it in terms of
local \v Cech data. Let $( \Omega^{(0,\bullet)}(X,E),\bar\partial)$
be the Dolbeault complex with values in the holomorphic vector
bundle $E$. We fix hermitian metrics on $T_X$ and on $E$. These
metrics induce an $L^2$ inner product $\langle\ ,\ \rangle$ on the
Dolbeault complex and a self-adjoint positive semidefinite Laplace
operator
$\Delta_{\bar\partial}=\bar\partial\bar\partial^*+\bar\partial^*\bar\partial$,
with discrete spectrum. By Hodge theory, the cohomology group
$H^j(X,E)$ is isomorphic to the space of harmonic forms
$\mathrm{Ker}(\Delta_{\bar\partial})$. Moreover we have the
following standard fact.

\begin{proposition} For any $D\in\cD_E(X)$ and $t>0$, $De^{-t\Delta_{\bp}}$
is a trace class operator on the Hilbert space of square integrable
Dolbeault forms. The expression
\[ \sum_{j=0}^n(-1)^j
\mathrm{tr}_{\Omega^{(0,j)}(X,E)}(D e^{-t\Delta_{\bar\partial}})
\] is independent of $t$ and is equal to $T_1(D)=L(D)$.
\end{proposition}

\begin{proof} We refer, e.g., to \cite{BGV} for the trace class property.
The independence of $t$ is checked by differentiation:
\begin{eqnarray*}
\lefteqn{\frac d{dt}\, \mathrm{tr}_{\Omega^{(0,j)}(X,E)}(D
e^{-t\Delta_{\bp}})= -\mathrm{tr}_{\Omega^{(0,j)}(X,E)}(D
e^{-t\Delta_{\bp}}(\bp\bp^*+\bp^*\bp))}&&
\\
&=& -\mathrm{tr}_{\Omega^{(0,j)}(X,E)}(D e^{-t\Delta_{\bp}}\bp\bp^*)
-\mathrm{tr}_{\Omega^{(0,{j-1})}(X,E)}(\bp D
e^{-t\Delta_{\bp}}\bp^*)
\\
&=& -\mathrm{tr}_{\Omega^{(0,j)}(X,E)}(D e^{-t\Delta_{\bp}}\bp\bp^*)
-\mathrm{tr}_{\Omega^{(0,{j-1})}(X,E)}( D
e^{-t\Delta_{\bp}}\bp\bp^*).
\end{eqnarray*}
Here we use the fact that $\bp$ commutes both with $D$ (since $D$ is
holomorphic) and with the Laplacian. Taking the sum with alternating
signs yields the claim.

Thus we can evaluate the sum in the limit $t\to\infty$. Since 0 is
an isolated eigenvalue of the positive semidefinite operator
$\Delta_{\bp}$ we obtain the alternating sum of traces on harmonic
forms, namely $L(D)$.
\end{proof}

\subsection{The $\sigma$-cocycle}
We introduce our main technical tool, a cocycle in a double complex
associated to an open set $U$. Here we describe its properties and
postpone its construction by heat kernel methods to Section
\ref{s-Valentine}.

Let $U$ be a sufficiently small open neighbourhood of an arbitrary
point in $X$. Let $A=\cD_E(U)$ and let $B=C^\infty(U)$ be the
algebra of smooth complex valued functions on $U$. Consider further
$C_p(A)=A\otimes \bar A^{\otimes p}$ with Hochschild differential
$b$ of degree $-1$ and $C_p(B)=B\otimes \bar B^{\otimes p}$ with
differential $s$ of degree $+1$ given by
\[
s(\rho_0\otimes\cdots\otimes\rho_p)=1\otimes\rho_0\otimes\cdots\otimes\rho_p.
\]
Let $C^c_p(B)$ be the subcomplex spanned by $(\rho_0,\dots,\rho_p)$
with compact common support $\cap_i\mathrm{supp}(\rho_i)$. Let us
denote by $[f(t)]_-=a_{-N}t^{-N}+\cdots+a_{-1}t^{-1}+a_0$ the
non-positive part of an asymptotic Laurent series $f(t)\sim
\sum_{j\geq -N}a_jt^{j}$ in $t$.
\begin{proposition}\label{p-xx} Let $U\subset X$ be an open set.
Let $A=\cD_E(U)$, $B=C^\infty(U)$. For each choice of hermitian
metrics on $T_X$ and $E$ there exist linear maps
\[
\sigma_p\colon C_p(A)\otimes C^c_p(B)\to\mathbb C[t^{-1}],
\]
such that the coefficients of $\sigma_p(D_0,\dots
D_p;\rho_0,\dots,\rho_p)$ are continuous in $(D_0,\dots,D_p)$ and
\begin{enumerate}
\item[(i)] Let $C_p^\varnothing(B)$ is the subcomplex spanned by
$(\rho_0,\dots,\rho_p)$ with empty common support
$\cap_{i=0}^{p}\mathrm{supp}(\rho_i)$. Then $\sigma_p$ vanishes on
$C_p(A)\otimes C_p^\varnothing(B)$.
\item[(ii)] For any $D\in C_{p+1}(A)$ and $\rho\in
C^c_p(B)$,
\[
\sigma_{p}(b D\otimes\rho) =\sigma_{p+1}(D\otimes s\rho),\qquad
p\geq 0,
\]
\item[(iii)]
$
\sigma_{0}(D,\rho)=\left[\sum_{j=0}^{n}(-1)^j\mathrm{tr}_{\Omega^{(0,j)}(U,E)}(\rho
De^{-t\Delta_{\bp}})\right]_-, $ $(n=\mathrm{dim}_{\mathbb C}(X))$.
\item[(iv)] Suppose that $U$ is some coordinate neighbourhood of a
point and let $c_E(U)$ be the cocycle appearing in Theorem
\ref{l-Miss Kenton}. Assume further that $\rho_0,\dots,\rho_{2n}\in
C^\infty_c(U)$ are functions such that the metrics on $T_X$ and $E$
are flat on some neighbourhood of  $\cap_{i=0}^{2n}
\mathrm{supp}(\rho_i)$. Then
\[ 
 \sigma_{2n}(c_E(U);\rho_0,\dots,\rho_{2n})=\frac{r 
 }{(2\pi
 i)^n}\int_U\rho_0d\rho_1\cdots d\rho_{2n},
\]
where $r$ is the rank of $E$.
\end{enumerate}
\end{proposition}

The proof is contained in Section \ref{s-Valentine}. We first show
how to use it to prove that $T_1$ is proportional to $T_3$.

\subsection{A local formula for the Lefschetz number}\label{ss-Nagasawa}
Here it is useful to replace the open cover considered in Section
\ref{s-Midori} by a refinement obtained from a triangulation of $X$.
Then the Hochschild--\v Cech cochains $(D^{(j)})$, constructed in
Section \ref{s-Midori} out of a global differential operator $D$,
define cochains, still denoted by $(D^{(j)})$ for the refinement.
Choose a smooth finite triangulation $|K|\to X$ of $X$, with
underlying simplicial complex $K$, with fixed total ordering of its
vertices. The {\em open star} of the triangulation is the open cover
$\mathcal U=(U_\alpha)_{\alpha\in K_0}$ of $X$ labeled by the set of
vertices of the triangulation, such that $U_\alpha$ is the
complement of the simplices not containing $\alpha$. By
construction, for all $\alpha_0<\cdots<\alpha_p$,
\begin{enumerate}
\item[(a)] $U_{\alpha_0}\cap\cdots \cap U_{\alpha_p}$ is empty or
contractible
\item[(b)] If $p>2n$ then $U_{\alpha_0}\cap\cdots \cap U_{\alpha_p}$ is empty.
\end{enumerate}
\begin{lemma}\label{l-Gil}
Let $(\rho_\alpha)$ be a partition of unity subordinate to the open
covering $(U_\alpha)$.
 Let $D\in\cD_E(X)$ and $s\in \check C^{2n}(\mathcal U,\mathcal
C_{2n}(\cD_E))$ be the cocycle \eqref{e-Toru Watanabe}. Then
\[
\sum_{p=0}^{2n}(-1)^{p}\mathrm{tr}(H^p(D))
=\sum_{\alpha_0<\dots<\alpha_{2n}}
\sigma_{2n}(s_{\alpha_0,\dots,\alpha_{2n}};\rho_{\alpha_0,\dots,\alpha_{2n}}).
\]
Here we use the abbreviation
\[
\rho_{\alpha_0,\dots,\alpha_q}=\sum_{\pi\in S_{q+1}}\sgn(\pi)
\rho_{\alpha_{\pi(0)}}\otimes\cdots\otimes\rho_{\alpha_{\pi(q)}}
\]
\end{lemma}

\begin{proof}
Out of $D$ we construct the cochains $D^{(j)}$ obeying
\eqref{e-staircase}.
\begin{eqnarray*}
T_1(D)&=&
\sum_{j=0}^n(-1)^j\left[\mathrm{tr}_{\Omega^{(0,j)}(X,E)}(D
e^{-t\Delta_\bp})\right]_-
\\
&=&\sum_\alpha
\sum_{j=0}^{n}(-1)^j\left[\mathrm{tr}_{\Omega^{(0,j)}(X,E)}(\rho_\alpha
D e^{-t\Delta_\bp})\right]_-
\\
&=&\sum_\alpha\sigma_0(D_\alpha;\rho_\alpha),\qquad
D_\alpha=D|_{U_\alpha}\in\cD_E(U_\alpha).
\end{eqnarray*}
Now $D_\alpha=b D^{(1)}_\alpha$ and Proposition \ref{p-xx} (ii)
implies
\begin{eqnarray*}
T_1(D)&=&\sum_\alpha\sigma_1(D^{(1)}_\alpha;1,\rho_\alpha)
\\
&=&\sum_{\alpha,\beta}\sigma_1(D^{(1)}_\alpha;\rho_\beta,\rho_\alpha)
\\
&=&\sum_{\alpha\neq\beta}\sigma_1(D^{(1)}_\alpha;\rho_\beta,\rho_\alpha)
+\sum_{\beta}\sigma_1(D^{(1)}_\beta;\rho_\beta,\rho_\beta)
\\
&=&
\sum_{\alpha\neq\beta}\sigma_1(D^{(1)}_\alpha-D^{(1)}_\beta;\rho_\beta,\rho_\alpha).
\end{eqnarray*}
In the last step we have replaced the last occurrence of
$\rho_\beta$ by $-\sum_{\alpha\neq\beta}\rho_\alpha \mod \mathbb
C1$. We see that $(\check\delta D^{(1)})_{\beta,\alpha}$ appears.
Thus we can iterate the procedure. At the $q$-th step we obtain
similarly for $q<2n$,
\begin{eqnarray*}
\sum_{\alpha_0<\cdots<\alpha_{q}}\sigma_q(\check \delta
D^{(q)}_{\alpha_0,\dots,\alpha_q};\rho_{\alpha_0,\dots,\alpha_q})
&=& \sum_{\alpha_0<\cdots<\alpha_{q}}\sigma_q(b
D^{(q+1)}_{\alpha_0,\dots,\alpha_q};\rho_{\alpha_0,\dots,\alpha_q})
\\
&=& \sum_{\alpha_0<\cdots<\alpha_{q}}\sigma_{q+1}(
D^{(q+1)}_{\alpha_0,\dots,\alpha_q};1\otimes\rho_{\alpha_0,\dots,\alpha_q})
\\
&=&\sum_{\alpha_0<\cdots<\alpha_{q+1}}\sigma_{q+1}(\check \delta
D^{(q+1)}_{\alpha_0,\dots,\alpha_{q+1}};\rho_{\alpha_0,\dots,\alpha_{q+1}})
\end{eqnarray*}
If $q=2n$ we have an additional term containing $s$ and we obtain
\begin{eqnarray*}
T_1(D)&=&\sum_{\alpha_0<\cdots<\alpha_{2n}}
\sigma_{2n}(s_{\alpha_0,\dots,\alpha_{2n}};\rho_{\alpha_0,\dots,\alpha_{2n}})
\\
&&+\sum_{\alpha_0<\cdots<\alpha_{2n+1}}\sigma_{2n+1}(\check \delta
D^{(2n+1)}_{\alpha_0,\dots,\alpha_{2n+1}};\rho_{\alpha_0,\dots,\alpha_{2n+1}}).
\end{eqnarray*}
Since there are no non-empty $(2n\!+\!2)$-fold intersections,
$(\rho_{\alpha_0},\dots,\rho_{2n+1})$ belongs to $C^\varnothing(B)$
and therefore, by Proposition \ref{p-xx}, (i), the second term
vanishes.
\end{proof}

Let us now choose the hermitian metrics so that they are flat on the
disjoint closed sets
$\cap_{j=0}^{2n}\mathrm{supp}(\rho_{\alpha_i})$,
$\alpha_0<\cdots<\alpha_{2n}$. By Proposition \ref{p-xx}, (iv), we
then obtain
\begin{eqnarray*} 
\lefteqn{
 \sum_{p=0}^{2n}(-1)^{p}\mathrm{tr}(H^p(D))}\\
 &=&(2n+1)!\frac {r 
 }{(2\pi i)^n}\sum_{\alpha_0<\dots<\alpha_{2n}}\lambda_{\alpha_0,\dots,\alpha_{2n}}(D)
 \int_X\rho_{\alpha_0}d\rho_{\alpha_1}\cdots d\rho_{\alpha_{2n}}.
\end{eqnarray*}
Now the common support of the functions $\rho_{\alpha_i}$ in each
summand is contained in a simplex
$\sigma_{\alpha_0,\dots,\alpha_{2n}}$. Moreover each of the
functions vanishes on the corresponding face and
$\sum_{i=0}^{2n}\rho_{\alpha_i}=1$ on some neighbourhood of the
simplex. Therefore the integral may be evaluated as follows.

\begin{lemma} Let $H_p\in\mathbb R^{p+1}$ be the hyperplane
$\sum_{i=0}^pt_i=1$ and $\Delta_p=H_p\cap [0,1]^{p+1}$ the standard
simplex, with (standard) orientation given by the ordered basis
$\partial_{t_1},\dots,\partial_{t_p}$. Let $\rho_0,\dots,\rho_p$ be
smooth functions defined on some open neighbourhood $U\subset H_p$
of $\Delta_p$ such that $\rho_0+\cdots+\rho_p=1$ and $\rho_i(t)=0$
if $t_i\leq0$. Then
\[
\int_{\Delta_p}\rho_0d\rho_1 \cdots d\rho_p=\frac1{(p+1)!}.
\]
\end{lemma}
\begin{proof}
We prove by induction in $p$ the more general formula
\[
\int_{\Delta_p}\rho_0^kd\rho_1\cdots
d\rho_p=\frac{k!}{(p+k)!},\qquad k=0,1,2,\dots
\]
This formula trivially holds for $p=0$. By the Stokes theorem,
\begin{eqnarray*}
\int_{\Delta_{p}}\rho_0^kd\rho_1\cdots d\rho_{p}&=&
-\int_{\Delta_{p}}\rho_0^kd\rho_1\cdots d\rho_{p-1}d\rho_{0}
\\
&=&
(-1)^{p}\frac1{k+1}\int_{\Delta_{p}}d(\rho_{0}^{k+1}d\rho_1\cdots
d\rho_{p-1})
\\
&=& (-1)^{p}\frac1{k+1}\int_{\partial\Delta_{p}}\rho_{0}^{k+1}
d\rho_1\cdots d\rho_{p-1}
\end{eqnarray*}
Since $\rho_j$ vanishes on the $j$th face of $\Delta_p$, only the
$p$th face (where $t_p=0$) contributes. This face is $\Delta_{p-1}$
and the restriction of $\rho_0,\dots,\rho_{p-1}$ obey the
assumptions of the lemma. Taking into account the sign $(-1)^p$
relating the orientation of $\Delta_{p-1}$ to the induced
orientation, we obtain
\[
\int_{\Delta_{p}}\rho_0^kd\rho_1\cdots d\rho_{p}=
\frac1{k+1}\int_{\Delta_{p-1}}\rho_{0}^{k+1} d\rho_1\cdots
d\rho_{p-1},
\]
proving the induction step.
\end{proof}

\begin{corollary}\label{c-Olivia}
Let $\epsilon(\alpha_0,\dots,\alpha_{2n})\in\{-1,1\}$ be the
orientation of the simplex $\sigma_{\alpha_0,\dots,\alpha_{2n}}$
relative to the canonical orientation of $X$. Then
\[ 
 T_1(D)=\frac {r 
 }{(2\pi i)^n}
 \sum_{\alpha_0<\dots<\alpha_{2n}}
\lambda_{\alpha_0,\dots,\alpha_{2n}}(D)\epsilon(\alpha_0,\dots,\alpha_{2n})
\]
\end{corollary}

\section{\dots and so is the second}
Let $T_2(D)=\int_X\chi_0(D)$ be the second trace. Let $C$ be the
cell decomposition of $X$ dual to the triangulation of subsection
\ref{ss-Nagasawa}. Its cells are in one-to-one correspondence with
the simplices of the triangulation. We denote
$C_{\alpha_0,\dots,\alpha_p}$ the $(2n-p)$-cell corresponding to the
simplex $\sigma_{\alpha_0,\dots,\alpha_{p}}$ with vertices
$\alpha_0,\dots,\alpha_p$. We orient the dual cells by the condition
that
$C_{\alpha_0,\dots,\alpha_p}\cdot\sigma_{\alpha_0,\dots,\alpha_p}=1$
on the intersection index (see Appendix \ref{a-Euen}).
\begin{proposition}
Let $s=s(D)$ be the cocycle  \eqref{e-Toru Watanabe}. Then
\[ 
 T_2(D)=
 \sum_{\alpha_0<\cdots<\alpha_{2n}}\int_{C_{\alpha_0,\dots,\alpha_{2n}}}\chi_{2n}(
 s_{\alpha_0,\dots,\alpha_{2n}}),
\]
where $\chi_{2n}$ is defined in Proposition \ref{p-Asfar} for the open
set $U_{\alpha_0}\cap\cdots\cap U_{\alpha_{2n}}$.
\end{proposition}

\begin{proof} We first prove by induction that for all $p=0,\dots,2n-1$,
\begin{equation}\label{e-John Lennon} 
 T_2(D)=
 \sum_{\alpha_0<\cdots<\alpha_p}\int_{C_{\alpha_0,\dots,\alpha_{p}}}\chi_{p}(b
 D^{(p+1)}_{\alpha_0,\dots,\alpha_{p}}),
\end{equation}
and then deduce the claim by doing a further induction step. For
$p=0$ eq.~\eqref{e-John Lennon} follows from
\[
T_2(D)=\sum_\alpha\int_{C(\alpha)}\chi_{0}(D|_{U_\alpha}),
\]
and $D|_{U_\alpha}=bD^{(1)}_\alpha$. Assume that the claim is proved
up to some $p< 2n-1$. 
 Then,
by Proposition \ref{p-Naoko} and the Stokes theorem (the signs are
discussed in the appendix, see \eqref{e-signs}), we get
\begin{eqnarray*} 
 T_2(D)&=&
 \sum_{\alpha_0<\cdots<\alpha_{p}}\int_{C_{\alpha_0,\dots,\alpha_{{p}}}}\chi_{{p}}(b
 D^{({p}+1)}_{\alpha_0,\dots,\alpha_{{p}}})
 \\
 &=&
 (-1)^p 
 \sum_{\alpha_0<\cdots<\alpha_{p}}\int_{C_{\alpha_0,\dots,\alpha_{{p}}}}d
 \chi_{{p}+1}(D^{({p}+1)}_{\alpha_0,\dots,\alpha_{{p}}})
 \\
 &=&(-1)^p
 (-1)^p 
 \sum_{\beta,\alpha_0<\cdots<\alpha_{p}}\int_{C_{\beta,\alpha_0,\dots,\alpha_{{p}}}}
 \chi_{{p}+1}(D^{({p}+1)}_{\alpha_0,\dots,\alpha_{{p}}})
 \\
 &=&
 \sum_{\alpha_0<\cdots<\alpha_{{p}+1}}\int_{C_{\alpha_0,\dots,\alpha_{{p}+1}}}
 \chi_{{p}+1}((\check\delta
 D^{({p}+1)})_{\alpha_0,\dots,\alpha_{{p}+1}}).
\end{eqnarray*}
Since $\check\delta D^{({p}+1)}=bD^{({p}+2)}$ if ${p}<2n-1$ 
the induction step is complete.

Now we do this step once more for ${p}=2n-1$. The calculation is the
same but the conclusion is different since $\check\delta
D^{(2n)}=s+b D^{(2n+1)}$. We obtain
\[ 
 T_2(D)=
 \sum_{\alpha_0<\cdots<\alpha_{2n}}\int_{C(\alpha_0,\dots,\alpha_{2n})}
 \chi_{2n}((s+b D^{(2n+1)})_{\alpha_0,\dots,\alpha_{2n}}).
\]
Moreover, $\chi_{2n}$ coincides with $\tau_{2n}$ composed with the
Taylor expansion and thus is a cocycle, i.e., it vanishes on exact
chains such as $bD^{(2n+1)}$.

The integral over the $0$-dimensional cycle
$C_{\alpha_0,\dots,\alpha_{2n}}$ is the evaluation of the integrand
times the sign of the orientation, that is the sign
$\epsilon(\alpha_0,\dots,\alpha_{2n})$ of the orientation of
$\sigma_{\alpha_0,\dots,\alpha_{2n}}$ relative to the orientation of
$X$.
\end{proof}

\begin{corollary}\label{c-Viola}
\[ 
T_2(D)=
r\sum_{\alpha_0<\cdots<\alpha_{2n}}\lambda_{\alpha_0,\dots,\alpha_{2n}}(D)
\epsilon(\alpha_0,\dots,\alpha_{2n}).
\]
\end{corollary}

\subsection{Proof of Theorem \ref{t-1}}
Recall that $T_1(D)=L(D)$ and that $T_2(D)=\int_X\chi_0(D)$. Theorem
\ref{t-1} follows from Corollary \ref{c-Olivia} and Corollary
\ref{c-Viola}. The missing step is the proof of Proposition
\ref{p-xx}, which appears in the next section.

\section{Asymptotic topological quantum mechanics}\label{s-Valentine}
In this section we prove Proposition \ref{p-xx} and give in
particular the construction of $\sigma_p$. Roughly speaking,
$\sigma_p$ is the cup product of a cochain $\Psi$, constructed using
topological quantum mechanics and a cochain $Z$ taking care of the
partition of unity. The formula for $\Psi$ is a version of the JLO
cocycle \cite{JLO} and is a regularized version of a cocycle
appearing in ``topological quantum mechanics'' \cites{Lysov,FLS}. It
is constructed with heat kernel methods. Here we need only the
asymptotic behaviour of these objects as time (or inverse
temperature \cite{JLO}) tends to zero, which allows us to replace
the heat kernel by a better behaved parametrix with support in a
neighbourhood of the diagonal.

We work in the context of Section \ref{s-1=3} and fix in particular
hermitian metrics on the holomorphic vector bundles $T^{1,0}X$ and
$E$.
\subsection{A parametrix for the heat equation}
We summarize here what we need about the heat kernel and refer to
\cite{BGV} for more details and proofs.  The heat operator
$e^{-t\Delta_{\bar\partial}}$ is an integral operator with kernel
$k_t\in\oplus_p\Gamma(X\times X, E^{0,p}\boxtimes (E^{0,p})^*)$,
where $E^{0,p}=\wedge^p(T^{0,1}X)^*\otimes E$: for any smooth
section $\phi\in \Omega^{0,\bullet}(X,E)$,
\[
e^{-t\Delta_{\bar\partial}}\phi(z)=\int_Xk_t(z,z')\cdot\phi(z')|dz'|,\qquad
t>0,
\]
is the solution of the heat equation
$\partial_tu+\Delta_{\bar\partial}u=0$ with initial data $\phi$.
Here $|dz'|$ denotes the Riemannian volume form. Let $d(z,z')$
denote the geodesic distance between $z,z'\in X$. Then the heat
kernel has an asymptotic expansion as $t\to 0$,
\begin{equation}\label{e-formal}
k_t(z,z')\sim \frac 1{(\pi
t)^n}e^{-\frac{d(z,z')^2}t}\left(\Phi_0(z,z')+t\Phi_1(z,z')+t^2\Phi_2(z,z')+\cdots\right).
\end{equation}
The smooth kernels $\Phi_j(z,z')$ can be chosen to vanish except on
an arbitrary small neighbourhood $d(z,z')<\varepsilon$ of the diagonal.
The precise meaning of the expansion is that if $k^N_t$ is the
truncation of the series at the $N$th term and $\|\cdot\|_\ell$
denotes the $C^\ell$-norm on sections of the hermitian bundle
$E^{0,p}\boxtimes (E^{0,p})^*$ on $X\times X$, then for all
$\ell,j,\alpha\geq0$ and $N$ sufficiently large, depending on
$\ell,j,\alpha$,
\begin{equation}\label{e-Vallorbe}
\|\partial_t^j(k_t-k^N_t)\|_\ell=O(t^\alpha),\qquad
\|(\partial_t+\Delta_{\bar\partial})k^N_t\|_\ell=O(t^\alpha).
\end{equation}
 Also, with the same hypotheses, for any smooth section $\phi$,
\begin{equation}\label{e-Dijon}
\lim_{t\to 0^+} \| K_t^N\phi-\phi\|_\ell=0,
\end{equation}
where $K_t^N$ denotes the integral operator with kernel $k_t^N$.

\subsection{Hochschild cohomology}

Let $A$ be an associative algebra with unit and let $(M=\oplus
M^j,d_M)$ be a complex of $A$-bimodules such that $M^j=0$ for all
but finitely many $j$. Recall that the Hochschild cochain complex
$C^\bullet(A,M)$ with values in $M$ is the total complex of the
double complex
\[
C^{p,q}(A,M)=\mathrm{Hom}(A^{\otimes p},M^q)
\]
and differential $\delta=d_H+(-1)^pd_M\colon C^{p,q}\to
C^{p+1,q}\oplus C^{p,q+1}$ with
\begin{eqnarray*}
d_H\varphi(a_1,\dots,a_{p+1})&=&
a_1\varphi(a_2,\dots,a_{p+1})\\
&&+\sum_{l=1}^{p}(-1)^l\varphi(a_1,\dots,a_la_{l+1},\dots,a_{p+1})\\
&&+(-1)^{p+1}\varphi(a_1,\dots,a_{p})a_{p+1}.
\end{eqnarray*}
The complex of $A$-bimodules  dual to $M$ is $(M^*=\oplus
(M^*)^j,d_{M^*})$ with $(M^*)^j=(M^{-j})^*$,
$d_{M^*}\varphi=(-1)^j\varphi\circ d_M$ for $\varphi\in (M^j)^*$ and
action of $A$ defined by $a\cdot\varphi(x)=\varphi(xa),\varphi\cdot
a(x)=\varphi(ax), a\in A, x\in M$. With these definitions,
$C^\bullet(A,A^*)$ is the complex dual to the Hochschild chain
complex $C_\bullet(A)$.

With any homomorphism $\bullet\colon M_1\otimes_A M_2\to M_3$ of
complexes of $A$-bimodules is associated a chain map, the {\em cup
product} $\cup\colon C^{p,q}(A,M_1)\otimes C^{p',q'}(A,M_2)\to
C^{p+q,p'+q'}(A,M_3)$,
\[
\varphi\cup\psi(a_1,\dots,a_{p+q})=
(-1)^{qp'}\varphi(a_1,\dots,a_p)\bullet\psi(a_{p+1},\dots,a_{p+q}).
\]
We will use this construction in two special cases: (a)
$M_1=M_2=M_3=M$ is a differential graded algebra whose product
factors through $M\otimes_A M$ defining thus a map $\bullet\colon
M\otimes_A M\to M$. (b) $M_1=M$ is a complex of $A$-bimodules,
$M_2=M^*$, $M_3=A^*$ with zero differential and $\bullet\colon
M^*\otimes_A M\to A^*$ is the map $(\varphi,x)\mapsto
(y\mapsto\varphi(xy))$.

\subsection{A JLO-type cocycle in the Hochschild--Dolbeault double complex}
Let $U$ be an open subset of $X$ and $A=\cD_E(U)$ be the algebra of
differential operators on the restriction of $E$ to $U$. The
Dolbeault complex
$(M_c(U)=\Omega_c^{0,\bullet}(U)\otimes_{\cO_X(U)}\cD_E(U),\bar\partial\otimes
\mathrm{id})$ with compact support and values in $\cD_E$ is a
locally convex differential graded algebra and an $A$-bimodule. In
local coordinates it is the graded algebra generated by
$M_r(C^\infty_c(U))$ of degree $0$, $d\bar z_i$ of degree 1 and
$\partial_{z_i}$ of degree zero.  The algebra $M_c(U)$ is the
inductive limit over $j$ and $K$ of the locally convex subalgebras
$M_{K,j}$ of operators of order at most $j$ and with support on a
compact subset $K\subset U$. The space $M_{K,j}$ is the space of
sections $x\to D_x$ of some vector bundle on $U$ with support in
$K$, and has the topology defined by the system of seminorms given
by the $C^\ell$-norms, for all $\ell$.

\begin{proposition}\label{p-Lili}
Let $U\subset X$ be an open subset, $A=\cD_E(U)$ and $M_c=M_c(U)$ be
the Dolbeault complex with values in $A$ and compact support. Let
$k^N_s$ be a parametrix, with support in some small neighbourhood of
the diagonal, obtained by truncating the formal series
\eqref{e-formal} at the $N$-th term. Suppose that $D_0\in M_c^p$,
$D_1,\dots,D_p\in A$. Then, for any sufficiently large $N$,
\begin{eqnarray*}
\lefteqn{\Psi_p(D_0,\dots,D_p)}\\
&&=(-1)^{\frac{p(p+1)}2}\left[\int_{t\Delta_p} \mathrm{Str}
(D_0K^N_{s_0}[\bar\partial^*,D_1]K^N_{s_1}\cdots
[\bar\partial^*,D_p]K^N_{s_p})ds_1\cdots ds_p\right]_-,
\end{eqnarray*}
where $\mathrm{Str}$ denotes the alternating sum of traces over the
Hilbert space of square integrable sections of $\wedge
(T^{0,1}U)^*\otimes E|_U$, is independent of $N$ for large $N$ and
defines a continuous cocycle
\[
\Psi=\sum_p \Psi_p\in \oplus_{p=0}^n\mathrm{Hom}(M^p_c\otimes \bar
A^{\otimes p},\mathbb C)[t^{-1}]\simeq C^0(A,M_c^*)[t^{-1}].
\]
\end{proposition}
\begin{proof}
The alternating trace $\mathrm{Str}
(D_0K^N_{s_0}[\bar\partial^*,D_1]K^N_{s_1}\cdots
[\bar\partial^*,D_p]K^N_{s_p})$ is the integral $\int_X\alpha_p|dx|$
of some function $\alpha_p\in C^\infty(X\times \Delta_p)$ with
support in some neighbourhood of the support of $D_0$. This function
has the form
\[
\alpha_p(x,s)=\int_{X^{p}}\mathrm{str}\left(D_0\prod_{j=0}^p
[\bar\partial^*,D_j]k^N_{s_j}(x_j,x_{j+1})\right)\prod_{j=1}^p|dx_{j}|,
\quad x_0=x_{p+1}=x,
\]
where the differential operators $[\bar\partial^*, D_j]$ act with
derivatives with respect to $x_j$ (the product is the composition of
linear maps in the conventional order). The supertrace
$\mathrm{str}$ is the alternating sum of traces over the fibres
$\wedge^jT_x^{0,1}X^*\otimes E_x$ at $x\in U$. The integral is
actually over a small neighbourhood of $(x,\dots,x)\in X^p$. Since
$k^N_s(z,z')$ is a smooth kernel, $\alpha_p(x,s)$ is smooth for $s$
in the interior of the simplex $t\Delta_p$. It is also continuous on
its boundary for any fixed $t$, uniformly in $x$, as can be seen
using \eqref{e-Dijon}. By rescaling $s=ts'$ we see that
$\int_{t\Delta_p}\alpha_p(x,s)\prod ds_i$ has an asymptotic
expansion as a Laurent series in $t$ whose singular part is not
affected by corrections of order $s^{N+1}$ to $k^N_s$ for large
enough $N$. Thus the expression for $\Psi_p$ is independent of $N$
for $N$ large enough. For further details see appendix \ref{a-Zhimei}.

The proof of the cocycle relation is similar to the proof in
\cite{JLO}. The Hochschild differential $d_H\Psi_p$ can be written
as the alternating sum of integrals of a differential form on
$X\times t\partial_i\Delta_{p+1}$, where $\partial_i\Delta_{p+1}$ is
the $i$-th face $s_i=0$ of the simplex $\Delta_{p+1}$. Using the
Stokes theorem and heat equation for $k^N_{s_i}$ (which holds up to
terms we can neglect by \eqref{e-Vallorbe}) to compute the
differential with respect to $s$ we obtain
\begin{gather*}
\Psi_p(D_0D_1,\dots,D_{p+1})-\Psi_p(D_0,D_1D_2,\dots,D_{p+1})+\cdots\\
+(-1)^{p+1}\Psi_p(D_{p+1}D_0,\dots,D_p)
=\Psi_{p+1}([\bar\partial,D_0],\dots,D_{p+1}),
\end{gather*}
which is the claim.
\end{proof}

\subsection{Construction of $\sigma_p$}
Let now $\rho_0,\dots,\rho_p\in C^\infty(U)$. View $C^\infty(U)$ as
a subalgebra of $M=\Omega^{0,\bullet}(U)\otimes_{\cO_X(U)}\cD_E(U)$
embedded as $C^\infty(U)\otimes\mathrm{id}$. Since $C^0(A,M)=M$, we
may consider $\rho_i$ as a 0-cochain and define
\[
Z^p(\rho_0,\dots,\rho_p)=
\rho_0\cup\delta\rho_1\cup\cdots\cup\delta\rho_p\in
C^p(A,M),
\]
where the cup product is defined using the product $M\otimes_A M\to
M$. Clearly
\begin{equation}\label{e-Z}
\delta Z^p(\rho_0,\dots,\rho_p)=Z^{p+1}(1,\rho_0,\dots,\rho_p).
\end{equation}
If $\cap_i\mathrm{supp}(\rho_i)$ is compact, then
$Z^p(\rho_0,\dots,\rho_p)$ takes values in differential operators
with compact support and therefore is a cochain in $C^p(A,M_c)$.

Let $\cup\colon C^\bullet(A,M_c^*)\otimes C^\bullet(A,M_c)\to
C^\bullet(A,A^*)$ be the cup product associated with the map
$M_c^*\otimes_A M_c\to A^*$ sending $\varphi\otimes x$ to the linear
form $a\mapsto\varphi(xa)$. We set
\[
\sigma_p(\rho_0,\dots,\rho_p)=\Psi\cup Z^p(\rho_0,\dots,\rho_p)\in
C^p(A,A^*)[t^{-1}].
\]
\subsection{Proof of Proposition \ref{p-xx}}
Claim (ii) follows from the fact that $\Psi$ is a cocycle and
eq.~\eqref{e-Z}. To prove the remaining claims let us write
$\sigma_p$ more explicitly:
\begin{eqnarray}\label{e-Maria}
 \lefteqn{ \sigma_p(D_0,\dots,D_p;\rho_0,\dots,\rho_p)}\\
 &=&
  \sum_{j=0}^p(-1)^{j(p-j)}
  \Psi_j(Z_{p-j}^p(D_{j+1},\dots,D_p;\rho_0,\dots,\rho_p)D_0,D_1,\dots,D_j).
\notag
\end{eqnarray}
The component $Z_{p-j}^p$ in $\mathrm{Hom}(\bar
A^{\otimes{p-j}},M_c^j)$ of $Z^p$ is given by
\begin{eqnarray*}
Z_{p-j}^p(D_{j+1},\dots,D_p;\rho_0,\dots,\rho_p)= \sum_{\pi\in
S_{p-j,j}}\sgn(\pi)\rho_0 B_{\pi(1)}(\rho_1)\cdots
B_{\pi(p)}(\rho_p),
\end{eqnarray*}
where $B_i(\rho)=[D_{j+i},\rho]$ for $i=1,\dots,p-j$, and
$B_i(\rho)=[\bar\partial,\rho]$ for $i=p-j+1,\dots,p$. From these
expressions it is clear that (i) and (iii) hold. For (iii) see also Appendix \ref{a-Zhimei}, Remark \ref{rem:jlo:replacebyfullhk_0}.
\dontprint{
\begin{lemma}\label{lem_Z}
The $Z_j^p$ satisfy the identity:
\begin{align*}
  &(-1)^j D_{j+1}Z_j^p(D_{j+2},\dots,D_{p+1};\rho_1,\dots,\rho_{p+1})
   \ +(-1)^{j+1}Z_j^p(\bar b(D_{j+1},\dots,D_{p+1});\rho_1,\dots,\rho_{p+1})\\
  &\ +(-1)^{p+1}Z_j^p(D_{j+1},\dots,D_p;\rho_1,\dots,\rho_{p+1})D_{p+1}
   =[\bar\partial,Z_{j-1}^p(D_{j+1},\dots,D_{p+1};\rho_1,\dots,\rho_{p+1})]\\
  &\ + Z_j^{p+1}(D_{j+1},\dots,D_{p+1};1,\rho_1,\dots,\rho_{p+1})
\end{align*}
\end{lemma}
} Let us turn to (iv). We need to evaluate
$\sigma_{2n}(c_E(U);\rho_0,\dots,\rho_{2n})$. By multiplying
$\rho_0$ by a partition of unity we may assume that the support of
$\rho_0$ is contained in a small coordinate neighbourhood of a
point. We have to compute a sum of $(2n)!$ terms of the form
\eqref{e-Maria} where $D_0=1$ and the remaining $D_k$ are partial
derivatives $\partial_{z_i}$ or operators of multiplication by
$z_i$. The arguments $D_k$ occurring in $Z_{2n-j}^{2n}$ appear in
the combination $[D_k,\rho_l]$ which vanishes if $D_k=z_i$.
Therefore the only non-vanishing terms in the sum \eqref{e-Maria}
have $j\geq n$ and $D_{j+1},\dots, D_{2n}$ are all derivatives
$\partial_{z_i}$. On the other hand, if $j>n$ then $Z^{2n}_{2n-j}$
vanishes since a product of more than $n$ $(0,1)$-forms is zero.
Thus only the term with $j=n$ survives and we have (setting
$\partial_i=\partial_{z_i}$)
\[
Z_n^{2n}(\partial_{i_1},\dots,\partial_{i_n};\rho_0,\dots,\rho_{2n})=
\rho_0\frac{\partial\rho_1}{\partial
z_{i_1}}\cdots\frac{\partial\rho_n}{\partial
z_{i_n}}\bar\partial\rho_{n+1}\cdots\bar\partial\rho_{2n}+\cdots
\]
where the dots denote the remaining shuffles. Therefore
\begin{equation}\label{eq-FatherMartin} 
 \sigma_{2n}(c_E(U);\rho_0,\dots,\rho_{2n})=
 (-1)^{n(n+1)/2}(-1)^n\sum_{\pi\in
 S_n}\sgn(\pi)\Psi_n(B,z_{\pi(1)},\dots,z_{\pi(n)}),
\end{equation}
where $B$ is the multiplication operator
\[
B=\sum_{\pi\in
S_{2n}}\sgn(\pi)\rho_0\frac{\partial\rho_{\pi(1)}}{\partial
z_{1}}\cdots\frac{\partial\rho_{\pi(n)}}{\partial
z_{n}}\frac{\partial\rho_{\pi(n+1)}}{\partial \bar
z_{1}}\cdots\frac{\partial\rho_{\pi(2n)}}{\partial \bar z_{n}}\,
d\bar z_1\wedge\cdots\wedge d\bar z_n.
\]
The sign $(-1)^{n(n+1)/2}$ is the sign of the permutation mapping
$(\partial_1,z_1,\dots,\partial_n, z_n)$ to
$(z_1,\dots,z_n,\partial_1,\dots,\partial_n)$; the sign $(-1)^n$ is
the sign appearing in \eqref{e-Maria} for $j=n$, $p=2n$. Note that
since $B$ is the operator of multiplication by a $(0,n)$-form, the
only trace appearing in the alternating sum defining $\Psi_{n}$ is
the trace over $\Omega^{0,n}$ and it comes with a sign $(-1)^n$. Let
us calculate $\Psi_n(B,z_1,\dots,z_n)$. The calculation for all
other permutations is similar and gives the same contribution to the
sum over $S_n$.
\begin{eqnarray*}\label{e-Befana}
\lefteqn{ \Psi_n(B,z_1,\dots,z_n)}\\
&=&(-1)^{n}(-1)^{n(n+1)/2}
\int_{t\Delta_n}\mathrm{tr}_{\Omega^{0,n}}(B K_{s_0}
[\bar\partial^*,z_1] K_{s_1} \cdots [\bar\partial^*,z_n]
K_{s_n})ds_1\cdots ds_n.
\end{eqnarray*}
With our assumption on the metrics, the heat kernel is the standard
heat kernel on $\mathbb C^n$. In this case,
\[
\bar\partial=\sum d\bar z_i\frac{\partial}{\partial \bar z_i},\qquad
\bar\partial^*=-\sum \frac{\partial}{\partial
z_i}\,\iota_{\scriptstyle{\frac\partial{\partial\bar
z_i}}},\quad,\Delta_{\bar\partial}=-\sum_{j=1}^n\frac{\partial^2}{\partial
z_j\partial\bar z_j},
\]
where $\iota$ denotes interior multiplication. Thus
$\Delta_{\bar\partial}$ is $-4$ times the standard Laplacian and the
kernel of $K_t$ is
\[
k_t(z,z')=\frac1{(\pi t)^n}e^{-\frac{|z-z'|^2}t}.
\]
Now $[\bar\partial^*,z_i]=-\iota_{\partial/\partial\bar z_i}$, which
commutes with $K_t$. The heat operators combine to $K_{s_0}\cdots
K_{s_n}=K_t$, since $\sum s_i=t$ on $t\Delta_n$. 
 The product $(-\iota_{\partial/\partial\bar z_1})\cdots (-\iota_{\partial/\partial\bar z_n})$
 acting on the basis $d\bar z_1\wedge\cdots\wedge d\bar z_n$ gives
 $(-1)^n(-1)^{n(n-1)/2}$.
Let us write $B=b(z)d\bar z_1\wedge\cdots\wedge d\bar z_n$. Then we
obtain
\begin{eqnarray*} 
 \Psi_n(B,z_{1},\dots,z_{n})&=&(-1)^n
 \int_U b(z)\mathrm{tr}_{\mathbb
 \C^r}k_t(z,z)|dz|\int_{t\Delta_{n}}ds_1\cdots ds_n
\\
 &=&\frac {(-1)^n
 r}{n!\pi^n}\int_Ub(z)|dz|.
\end{eqnarray*}
The standard volume form $|dz|$ is
\begin{eqnarray*}
|dz|&=&(-2i)^{-n}dz_1\wedge d\bar z_1\wedge\cdots\wedge
dz_n\wedge d\bar z_n\\
&=&(-2i)^{-n}(-1)^{n(n-1)/2}dz_1\wedge\cdots\wedge dz_n\wedge d\bar
z_1\wedge\cdots\wedge d\bar z_n.
\end{eqnarray*}
Thus $b(z)|dz|=(2i)^{-n}(-1)^{n(n+1)/2}\rho_0d\rho_1\cdots
d\rho_{2n}$. Inserting this in the formula \eqref{eq-FatherMartin}
gives the formula that had to be proved.
\appendix

\section{Triangulations and signs}\label{a-Euen}
 Let $T$ be a smooth finite triangulation of
the oriented $d$-dimensional manifold $X$. Let
$\sigma_{\alpha_0,\dots,\alpha_p}\subset X$ denote the simplex with
vertices $\alpha_0,\dots,\alpha_p$. It is the image of the standard
oriented simplex $\Delta_p=\{t\in[0,1]^{p+1}\,|\,\sum t_i=1\}$
sending the $i$-th vertex with $t_i=1$ to $\alpha_i$ and thus comes
with an orientation, for which
\begin{equation}\label{e-Eddie}
\partial\sigma_{\alpha_0,\dots,\alpha_p}=
\sum_{j=0}^p(-1)^j\sigma_{\alpha_0,\dots,\hat\alpha_j,\dots,\alpha_p}.
\end{equation}
 The cells of the dual cell decomposition $T^*$ of $X$ (see
\cite{Lefschetz}) are in one-to-one correspondence with the
simplices of the triangulation. The $(d-p)$-cell
$C_{\alpha_0,\dots,\alpha_p}$ intersects only the $p$-simplex
$\sigma_{\alpha_0,\dots,\alpha_p}$ and meets it transversally in
exactly one interior point. Let us orient the cells by the condition
that the intersection index is one:
\begin{equation}\label{e-Cat}
C_{\alpha_0,\dots,\alpha_p}\cdot\sigma_{\alpha_0,\dots,\alpha_p}=1
\end{equation}
This means in particular that the top-dimensional cells $C_\alpha$
have the same orientation as $X$. With this convention both
$C_{\alpha_0,\dots,\alpha_p}$ and $\sigma_{\alpha_0,\dots,\alpha_p}$
change their orientation under permutation of the indices according
to the sign of the permutation.

If $c_p$ is a $p$-cell of $T^*$ and $c'_{d-p+1}$ is a
$(d\!-\!p\!+\!1)$-cell of $T$, we have
\[
\partial c_p\cdot c'_{d-p+1}=(-1)^pc_p\cdot\partial c'_{d-p+1}.
\]
By combining this equation with \eqref{e-Eddie} and \eqref{e-Cat} we
obtain the formula for the boundary of dual cells:
\begin{equation}\label{e-signs}
\partial
C_{\alpha_0,\dots,\alpha_p}=(-1)^{d+p}\sum_{\beta}C_{\beta,\alpha_0,\dots,\alpha_p},
\end{equation}
with summation over all $\beta$ such that
$\beta,\alpha_0,\dots,\alpha_p$ are the vertices of a simplex of the
triangulation.

\section{Heat kernel estimates and asymptotic expansion}\label{a-Zhimei}

In this section, we show the existence of the asymptotic expansion in the definition of the JLO-cocycle (see Proposition \ref{p-Lili}). In the first subsection, it is shown that the integrand in the definition of the JLO-cocyle is smooth for $s\in[0,1]^{p+1}\setminus\{0\}$. In the second subsection, we apply this result to compute the asymptotic expansion. In particular it will follow from this computation that $\Psi_p$ is well defined and continuous in the operators $D_0,\dots,D_p$.

\subsection{Heat kernel approximation}\label{SS:heat_approx}

In order to show that the integrand $f(s)$ in the formula for $\Psi_p$ is smooth for $s\in[0,1]^{p+1}\setminus\{0\}$, we need some estimates for the approximated heat kernel. We recall from \cite{BGV} the notions of a generalized Laplacian and the corresponding heat kernel. A {\em generalized Laplacian} $H$ acting on sections of a vector bundle $\E$ over a  $d$-dimensional Riemannian manifold $(X,g)$ is a second-order differential operator, which in local coordinates can be written as $H=-\sum_{i,j=1}^d g^{ij} \partial_i\partial_j +\text{ first order terms}$. It is easy to verify that $\Delta_{\barpart}= \barpart\barpart^*+ \barpart^*\barpart$ is 4 times such a Laplacian if we set $\E=E\otimes \Lambda^{\bullet}T^{*(0,1)}X$. Therefore we may directly use the results about the heat kernel from \cite{BGV} considering $X$ as a smooth $2n$-dimensional Riemannian manifold.

We write $\DO_\E(X)$ for the space of smooth differential operators acting on smooth sections $\Gamma(X,\E)=\oplus\Omega^{0,j}(X)$. $\Gamma(X,\E)$ is a locally convex space where the norms are the $C^k$-norms. These norms can be constructed by choosing a finite open cover of coordinate neighbourhoods of $X$. We then consider a cover of $X$ of compact sets that are slightly smaller than the previous open sets. The $C^k$-norms are then defined by the sum of the $C^k$-norms on the compact sets and with respect to the corresponding coordinates. Furthermore we can assume that the $C^k$-norms on $\Gamma(X,\E)$ are increasing, i.e. $\|\phi\|_k\le \|\phi\|_\ell$ for $k\le\ell$.

The spaces $\DO_E^j(X)\subset \DO_E(X)$ of differential operators of order $j$ are spaces of sections of a certain hermitian vector bundle over $X$, and so one can define increasing $C^k$-norms on them in a similar way as above. Then
$\DO_\E(X)$ is an LF-space which is the strict inductive limit of $\DO^j_\E(X)$, see, e.g., \cite{Treves}.

For two vector bundles $\E_1$ and $\E_2$ over the manifold $X$, we denote by $\E_1\boxtimes\E_2$ the external tensor product which is a vector bundle over $X\times X$. The {\em heat kernel} $k_t(x,y)$ is a family of sections $k_t\in\Gamma(X\times X,\E\boxtimes \E^*)$ defined for $t>0$ which is $C^1$ with respect to $t$ and $C^2$ with respect to $x$ and solves the equation
\[
  \partial_tk_t(x,y) + \Delta_{\barpart} k_t(x,y) = 0
\]
with initial condition $\lim\limits_{t\to 0}k_t(x,y) = \delta(x-y)$ where $\delta$ is the Dirac distribution with respect to the Riemannian density on $X$. The heat kernel exists and is unique. There is an approximation to the heat kernel $k^N_t(x,y)$ of the form
\[
  k^N_t(x,y) = (\pi t)^{-n}e^{-d(x,y)^2/t}\sum_{i=0}^N t^i\Psi_i(x,y)
\]
where $d(x,y)$ is the geodesic distance and $\Psi_i(x,y)$ are linear maps $\E_y\to\E_x$ depending smoothly on $(x,y)$ and with support in the set where $d(x,y)\le\varepsilon$ for some fixed $\varepsilon$ which can be chosen arbitrarily small. Furthermore $\Psi_0(x,x)$ is the identity and  the maps $\Psi_i$  can be chosen so that the following theorem holds:

\begin{thm}\label{thm:jlo:hkapprox}Let here $\|\;.\; \|_\ell$ be $C^\ell$-norms for sections in the bundle $\E\boxtimes \E^*$.
\begin{itemize}
\item[(i)] $k^N_t$ approximates the heat kernel $k_t$ in the sense that
\[
  \|\partial_t^m(k_t-k^N_t)\|_\ell = \mathcal{O}(t^{N-n-\ell/2-m})\quad\text{for }t\to0\,.
\]
\item[(ii)] $k^N_t$ is an approximate solution of the heat equation such that the remainder $r^N_t(x,y):=(\partial_t+\Delta_\barpart)k^N_t(x,y)$ satisfies the estimates
\[
  \|\partial^k_tr^N_t\|_\ell<Ct^{N-n/2-k-\ell/2}
\]
for some constant $C$ depending on $\ell$ and $k$.
\end{itemize}
\end{thm}
\begin{proof}
See \cite{BGV}, theorem 2.23 or 2.30 for part (i) and theorem 2.20 for part (ii).
\end{proof}

\begin{remark}\label{rem:jlo:replacebyfullhk_0}
For any $D_0\in M_c^0$ we have
\[
  \Psi_0(D_0):=[\Str(D_0K^N_t)]_-=[\Str(D_0K_t)]_-\;.
\]
This follows directly from the estimates about the approximated heat kernel in part i) of the above theorem. This remark will be generalized for $\Psi_p$, $p=0,\dots,2n$ in remark \ref{rem:jlo:replacebyfullhk}.
\end{remark}

We define the operator $K_t$ on smooth sections $\varphi\in\Gamma(X,\E)$ by
\begin{equation}\label{E:K_t}
  (K_t\varphi)(x) = \int_X k_t(x,y)\varphi(y)|dy|_g
\end{equation}
where $|dy|_g$ is the Riemannian density on X. $\varphi_t := K_t\varphi$ is a solution of the heat equation $\partial_t\varphi_t + \Delta_{\barpart}\varphi_t = 0$ with initial condition $\lim\limits_{t\to 0}\varphi_t = \varphi$. In the same way, we also define the operators $K^N_t$ that correspond to the approximated heat kernel $k^N_t$. We also set $K_0=K_0^N=\mathrm{Id}$.

The operator $K_t$ satisfies the following estimates:

\begin{lemma}\label{lem:jlo:estimates}
We write $\|.\|_\ell$, $\ell\ge0$ for the $C^\ell$-norms on $\Gamma(X,\E)$ or $\DO_\E(X)$. Fix a $\delta>0$ small enough. Then for each $\ell$ and each of the following inequalities there is a constant $C$ so that  for all $s,s'\in[\delta,1]$ and $t\in[0,1]$,
\begin{alignat*}{2}
  (i)& \quad\ &&\|K_t^N\varphi-\varphi\|_\ell \le C\,\|\varphi\|_{\ell+1}\,\sqrt{t} \hspace{4cm}\\
  (ii)&      &&\|K_s^N\varphi-K_{s'}^N\varphi\|_\ell \le C\,\|\varphi\|_0\,|s-s'| \\
  (iii)&     &&\|K_s^N\varphi\|_\ell \le C\,\|\varphi\|_0\\
  (iv)&      &&\|DK_0^N\varphi\|_\ell = \|D\varphi\|_\ell \le C\,\|D\|_\ell\,\|\varphi\|_{\ell+d}
\end{alignat*}
for every differential operator $D\in\DO_\E(X)$ of degree $d$.
\end{lemma}
\begin{proof}
(i) The proof is essentially the same as for the first part of Theorem 2.29 in \cite{BGV}. We consider the formula (\ref{E:K_t}) for $K^N_t\varphi$, change to exponential coordinates for $y$ ($y\mapsto\exp_x y$) and write $\varphi(x,y) := \varphi(\exp_x y)$ and $\Psi_j(x,y)=\Psi_j(x,\exp_x y)$, in the latter case with a slight abuse of the notation. We may assume that $\varepsilon$ in the definition of $k^N_t$ is smaller than the injectivity radius of the exponential map, so that the previous change to exponential coordinates in well defined. The substitution $y=\sqrt{t}v$ leads to
\[
  (K_t^N\varphi-\varphi)(x) = \frac{1}{\pi^n} \int\limits_{T_xX} e^{-\|v\|^2} \left( \sum_{j=0}^N t^j \Psi_j(x,\sqrt{t}v)\varphi(x,\sqrt{t}v)\rho(x,\sqrt{t}v)-\varphi(x,0)\right) dv
\]
where we used $\frac{1}{\pi^n}\int_{T_xX} e^{-\|v\|^2}=1$, and $\rho(x,y):=\sqrt{\det(\exp_x^*g(y))}$ is the factor coming from the Riemannian density. As $\Psi_j(x,y)=0$ for $\|y\|>\varepsilon$, it is a compactly supported function on $TX$. For $j>0$, it is therefore clear -- by taking the supremum over $y$ -- that $t^j\Psi_j(x,y)\varphi(x,y)\rho(x,y)$ is bounded by a constant times $\sqrt{t}\|\varphi\|_0$. For $j=0$, we write $f(x,y)=\Psi_0(x,y) \varphi(x,y)\rho(x,y)$. As $f(x,0)=\varphi(x,0)$, we get by the mean value theorem for the $t^0$-term
\[
  \frac{1}{\pi^n} \int_{T_xX} v e^{-\|v\|^2} \,\partial_y f(x,\sqrt{t'}v) \sqrt{t} dv
\]
for some $t'\in[0,t]$. This expression is bounded by a constant times $\sqrt{t}\|\varphi\|_1$ and the claim follows in the case $\ell=0$. For $\ell>0$, we use the same arguments, but the function $f(x,y)$ is replaced by $\partial^\alpha_xf(x,y)$ where $|\alpha|\le\ell$.

(ii) $K_s^N$ is an integral operator with kernel with $C^1$-dependence on $s$ for $s>0$. Therefore the mean value theorem tells us that
\[
  |\partial_x^\alpha k_s^N(x,y)-\partial_x^\alpha k_{s'}^N(x,y)| \le \sup\limits_{s\in[\delta,1]} |\partial_s\partial_x^\alpha k_s^N(x,y)|\,|s-s'|
\]
from which the claim follows.

(iii) Is obvious as the kernel is smooth in $x$ for all $s\in[\delta,1]$.

(iv) Also obvious.
\end{proof}

By iterating the above lemma, we find the following estimate:

\begin{lemma}\label{lem:jlo:abschaetzung}
Let $D_i\in\DO_\E(X)$ be differential operators of degree $d_i$, $i=1,\dots,m$. Fix a $\delta>0$ and a set $I\subset\{1,\dots,m\}$. Let $s_i\in [0,1],s'_i=0$ for $i\in I$ and $s_i, s'_i\in [\delta,1]$ for $i\notin I$. Then there is a constant $C$ and an $L\le \ell+m+\sum_{i=1}^m d_i$ so that
\begin{equation*}
  \begin{split}
  &\|D_1K_{s_1}^ND_2K_{s_2}^N\cdots D_mK_{s_m}^N\varphi - D_1K_{s̈́_1'}^ND_2K_{s_2'}^N\cdots D_mK_{s_m'}^N\varphi\|_\ell \\ &\le C \|\varphi\|_L \left( \sum_{i\in I}\sqrt{s_i} + \sum_{i\notin I}|s_i-s_i'|\right) \prod_{j=1}^k \|D_j\|_L
  \end{split}
\end{equation*}
\end{lemma}

\begin{proof}
Using the triangle inequality and Lemma \ref{lem:jlo:estimates}, we find
\begin{align*}
  \|DK_s^N\varphi_1-DK_{s'}^N\varphi_2\|_\ell &\le \|(DK_s^N-DK_{s'}^N)\varphi_1\|_\ell + \|DK_{s'}^N(\varphi_1-\varphi_2)\|_\ell\\
&\begin{cases}
  \overset{s'=0}\le C\sqrt{s}\|\varphi_1\|_{\ell+d+1} \|D\|_\ell + C\|\varphi_1-\varphi_2\|_{\ell+d}\|D\|_\ell\\
  \overset{s'\ge\delta}\le C|s-s'|\,\|D\|_\ell\|\varphi_1\|_0 + C\|\varphi_1-\varphi_2\|_0 \|D\|_\ell\,.
\end{cases}
\end{align*}
The proof is straightforward by induction on $m$.
\end{proof}

\begin{lemma}\label{lem:jlo:psi_welldef}
The function $f(s)$, which is the integrand in the definition of $\Psi_p$, is continuous for $s\in[0,1]^{p+1}\setminus\{0\}$. In particular the integral over $t\Delta_p$ in the definition of $\Psi_p$ (see Proposition \ref{p-Lili}) is well defined for $t\in(0,1]$.
\end{lemma}
\begin{proof}
An operator $D$ on $\Gamma(X,\E)$ with continuous kernel $D(x,y)\in\Gamma(X\times X,\E\boxtimes\E^*)$ is of trace class, and the supertrace can be written as
\begin{align*}
  \Str(D) &= \sum_{k=0}^n (-1)^k\Tr_{\Omega^{0,k}(X,E)}(D)\\
  &= \sum_{k=0}^n (-1)^k \int_X\tr_{E_x\otimes\Lambda^{0,k}(T_xX)}D(x,x)|dx|_g \,.
\end{align*}
For the integral over $t\Delta_p$ in the definition of $\Psi_p$ to be convergent, it is sufficient to show that the function $f(s):=\Str(D_0 K_{s_0}^N [\barpart^*,D_1]K_{s_1}^N \cdots [\barpart^*,D_p] K_{s_p}^N)$ is continuous in $s$ for $s\in t\Delta_p$. As the heat kernel $k_{s_i}^N$ is $C^1$ with respect to $s_i$ for $s_i>0$, this is clear except for points on the boundary of $t\Delta_p$. For a point $s'\in t\partial\Delta_p$, let $I$ be the subset of $\{1,\dots,n\}$ so that $s_i'=0 \Leftrightarrow i\in I$ and take a $\delta>0$ so that $s_i'>\delta$ $\forall i\notin I$. As there is at least one $i\notin I$ and as the trace is cyclic, we can w.l.o.g. assume that $p\notin I$. To simplify the notation, we set $A_{s_0\dots s_{p-1}} = D_0 K_{s_0}^N [\barpart^*,D_1]K_{s_1}^N \cdots [\barpart^*,D_{p-1}] K_{s_{p-1}}^N$ and $B_{s_p} = [\barpart^*,D_p] K_{s_p}^N$. We write the supertrace as
\[
  \Str(D_0 K_{s_0}^N\cdots [\barpart^*,D_p] K_{s_p}^N) = \!\!\sum_{\abtop{k=0\dots n}{i=1\dots i_k}} \!(-1)^k \!\!\! \int\limits_{X\times X}\!\! {\langle v^k_i|A_{s_0\dots s_{p-1}}(x,y) B_{s_p}(y,x) | v^k_i \rangle} \,dx\,dy
\]
where $\{v^k_i\}$ for fixed $k$ and $i=1,\dots, i_k$ is a basis for $E\otimes\Lambda^{0,k}(T_xX)$. Now we consider $A_{s_0\dots s_{p-1}}$ as operator acting on $\varphi_{s_p}:= B_{s_p}(\,\cdot\,,x)v$ where $x\in X$ and $v\in\E_x$ are considered as parameter. Then we get by the triangle inequality and Lemma \ref{lem:jlo:abschaetzung} that
\begin{equation*}\label{F:jlo:big'n'ugly}
\begin{aligned}
  &\|A_{s_0\dots s_{p-1}}\varphi_{s_p}-A_{s_0'\dots s_{p-1}'}\varphi_{s_p'}\|_0 \\
  &\le \| (A_{s_0\dots s_{p-1}}-A_{s_0'\dots s_{p-1}'})\varphi_{s_p} \|_0 + \| A_{s_0'\dots s_{p-1}'}(\varphi_{s_p}-\varphi_{s_p'}) \|_0\\
  &\le \widetilde{C}\left( \sum_{i\in I\setminus\{p\}} \sqrt{s_i} + \sum_{i\notin I\cup\{p\}} |s_i'-s_i|\right)\|\varphi_{s_p}\|_L  + \widetilde{C} \|\varphi_{s_p}-\varphi_{s_p'}
\end{aligned}
\end{equation*}
where $\widetilde{C} = C \|D_0\|_L\prod_{j=1}^{p-1} \|[\barpart^*,D_j]\|_L$.
We use the mean value theorem and find
\begin{align*}
  \|\varphi_s-\varphi_{s'}\|_L &\le |s-s'| \sup\limits_{\abtop{(x,v)\in\E, \|v\|\le 1}{u\in [\delta,1]}} \|B_u(\,\cdot\,,x)v\|_{L+1}\\
  &\le C |s-s'|\, \|[\barpart^*,D_p]\|_{L+1}\,.
\end{align*}
As the integral of the trace is over a compact set, we have shown that $f$
is continuous in $s$ for $s\in t\Delta$ and
\begin{equation}\label{F:jlo:bigestimate}
  |f(s)-f(s')| \le C \left( \sum_{i\in I} \sqrt{s_i} + \sum_{i\notin I} |s_i'-s_i|\right) \prod_{j=0}^p \|D_j\|_{L+2}
\end{equation}
\end{proof}

\begin{proposition}\label{prop:jlo:fsmooth}
The function $f(s)$ (see Lemma \ref{lem:jlo:psi_welldef}) is $k$-times continuously differentiable for $s\in[0,1]^{p+1}\setminus\{0\}$ and $N=N_k$ large enough.
\end{proposition}
\begin{proof}
The proof works in exactly the same way as in the previous lemmata (\ref{lem:jlo:estimates}, \ref{lem:jlo:abschaetzung}, \ref{lem:jlo:psi_welldef}). We  generalize the estimates in Lemma \ref{lem:jlo:estimates} by adding time derivatives: Fix a $\delta>0$ and assume $s,s'\in[\delta,1]$ and $t\in[0,1]$. Then for each $\ell,m$ and each of the following inequalities there is a constant $C$ so that
\begin{alignat*}{2}
  (i)  & \quad\ &&\|\partial_t^mK_t^N\varphi-(-\Delta_\barpart)^m\varphi\|_\ell \le C\,\|\varphi\|_{2m+\ell+1}\,\sqrt{t} \hspace{4cm}\\
  (ii) &      &&\|\partial_s^mK_s^N\varphi-\partial_{s'}^mK_{s'}^N\varphi\|_\ell \le C\,\|\varphi\|_0\,|s-s'| \\
  (iii)&      &&\|\partial_s^mK_s^N\varphi\|_\ell \le C\,\|\varphi\|_0\\
  (iv) &      &&\|(-\Delta_\barpart)^m D K_0^N\varphi\|_\ell = \|(-\Delta_\barpart)^m D\varphi\|_\ell \le C\,\|\varphi\|_{\ell+d+2m}
\end{alignat*}
where $\|.\|_\ell$, $\ell\ge0$ are $C^\ell$-norms on $\Gamma(X,\E)$, $\DO_\E(X)$ respectively. We only prove the first estimate as the others are easy to show (see Lemma \ref{lem:jlo:estimates}). Recall from Theorem \ref{thm:jlo:hkapprox} that the remainder $r_t^N=(\partial_t+\Delta_\barpart)k^N_t$ satisfies $\|\partial_t^kr_t^N\|_\ell<Ct^{N-k-(n+\ell)/2}$. By the iterated application of $\partial_tk^N_t=-\Delta_\barpart k^N_t + r_t^N$, we find
\[
  \partial_t^mk_t^N=(-\Delta_\barpart)^mk_t^N+\sum_{j=0}^{m-1}(-\Delta_\barpart)^{m-1-j}\partial_t^jr_t^N
\]
and hence the estimate
\begin{gather*}
  \|\partial_t^mK_t^N\varphi-(-\Delta_\barpart)^mK_t^N\varphi\|_\ell \le \sum_{j=0}^{m-1}\|\Delta_\barpart^{m-1-j}\partial_t^jr_t^N\varphi\|_\ell\\
  \le\sum_{j=0}^{m-1}\|\Delta_\barpart^{m-1-j}\|_\ell\|\partial_t^jr_t^N\|_{\ell+2(m-1-j)}\|\varphi\|_0 \le C\|\varphi\|_0t^{N-m+1-(n+\ell)/2}\,.
\end{gather*}
We require $N$ to be large enough, namely $N\ge\frac{n+\ell-1}{2}+m$. On the other hand we have the estimate
\[
  \|(-\Delta_\barpart)^mK_t^N\varphi-(-\Delta_\barpart)^m\varphi\|_\ell \le C\|\Delta_\barpart^m\|_\ell \|K_t^N\varphi-\varphi\|_{\ell+2m}\le C\|\varphi\|_{2m+\ell+1}\sqrt{t}
\]
The estimate $(i)$ then follows by the triangle inequality.

Using the above estimates, it is now straightforward to generalize Lemma \ref{lem:jlo:abschaetzung} to
\begin{equation*}
  \begin{split}
  &\|D_1\partial_{s_1}^{m_1}K_{s_1}^ND_2\partial_{s_2}^{m_2}K_{s_2}^N\cdots D_p\partial_{s_p}^{m_p}K_{s_p}^N\varphi - D_1\partial_{s_1'}^{m_1}K_{s_1'}^ND_2\partial_{s_2'}^{m_1}K_{s_2'}^N\cdots D_m\partial_{s_p'}^{m_p}K_{s_p'}^N\varphi\|_\ell \\ &\le C \|\varphi\|_L \left( \sum_{i\in I}\sqrt{s_i} + \sum_{i\notin I}|s_i-s_i'|\right)
  \end{split}
\end{equation*}
which is true for some $L\le\ell+\sum_i(d_i+2m_i)+1$. Then we see as in Lemma $\ref{lem:jlo:psi_welldef}$ that the partial derivatives of $f(s)$ up to degree $k$ are continuous.
\end{proof}

\subsection{Computation of $\Psi_p$ and power counting}

In this subsection, we explain an algorithm to compute $\Psi_p$ which will lead to the result summarized in the following proposition:

\begin{proposition}\label{prop:jlo:comp_power}
Let $n$ be the complex dimension of $X$. Take the operators $D_0, D_1,\dots,D_p$ as in Proposition \ref{p-Lili}. We write $d$ for the sum of the degrees of the differential operators $D_0, [\barpart^*,D_1],\dots,[\barpart^*,D_p]$ which are defined on a small (see remark below) open set $U\subset X$. Recall that the approximated heat kernel depends on the constants $N$ and $\varepsilon$. Then for $N$ big enough and $\varepsilon$ small enough, $\Psi_p(D_0,\dots,D_p)$ is well defined and a polynomial in $t^{-1}$ of degree $n-p+\lfloor\frac d2\rfloor$. More precisely, for $N\ge n-p+\lfloor\frac d2\rfloor$ and $\varepsilon<\frac{1}{p+1}\dist(X\setminus U,\supp(D_0))$, where $\dist$ means the geodesic distance, it is independent of $N$ and $\varepsilon>0$. Furthermore, $\Psi_p(D_0,\dots,D_p)$ depends continuously on $D_0,D_1,\dots,D_p$.
\end{proposition}

\begin{remark}
The set $U$ in the above proposition has to be small in the sense that Lemma \ref{lem:jlo:comp2} holds for any compact subset $K\subset U$. For larger $U$ the above proposition would still be true with the exception that the upper bound for $\varepsilon$ would need a more careful definition.
\end{remark}

The main idea of the computation is to ``move'' the operators
$[\barpart^*,D_i]$ in the formula for $\Psi_p$ (see Proposition
\ref{p-Lili}) to the left and to use a saddle point approximation
for the heat kernel integrals. As a preparation for this
computation, we formulate the following three lemmata:

\begin{lemma}\label{lem:jlo:comp1}
Let $U\subset X$ be a open subset of $X$ so that the exponential map w.r.t. any point in $U$ and restricted to the preimage of $U$ is a diffeomorphism. Assume that $x_1,x_2\in U$, then (in local coordinates) there is a smooth matrix valued function $a(x_1,x_2)$ so that
\[
  \frac{\partial}{\partial x_2} d(x_1,x_2)^2 = a(x_1,x_2)\frac{\partial}{\partial x_1} d(x_1,x_2)^2\,.
\]
\end{lemma}
\begin{proof}
We construct such a map explicitly: We introduce the coordinates $(x,\xi)=(x_1,\log_{x_1}x_2)$ and $(y,\eta)=(x_2,\log_{x_2}x_1)$. Obviously $|\xi|=d(x_1,x_2)=|\eta|$. Therefore we find
\begin{align*}
  \frac{\partial}{\partial x_2} d(x_1,x_2)^2 = \frac{\partial\xi}{\partial x_2}\frac{\partial}{\partial \xi} |\xi|^2 = \frac{\partial\xi}{\partial x_2}\frac{\partial{\eta}}{\partial\xi}\frac{\partial}{\partial\eta} |\eta|^2 = \frac{\partial\xi}{\partial x_2}\frac{\partial{\eta}}{\partial\xi}\frac{\partial x_1}{\partial\eta}\frac{\partial}{\partial x_1} d(x_1,x_2)^2\,.
\end{align*}
As the exponential coordinates are smooth coordinates, the lemma follows.
\end{proof}

\begin{lemma}\label{lem:jlo:comp2}
Let $K\subset X$ be a sufficiently small compact neighbourhood of
any point, so that the exponential map w.r.t. any point in $K$ and
restricted to the preimage of $K$ is injective. Take $x_1,x_2,x_3\in
K$ and $s_1,s_2\in (0,1]$, then for fixed $x_1,x_3,s_1,s_2$ the
function
\[
  f(x_2) = \frac{d(x_1,x_2)^2}{s_1}+\frac{d(x_2,x_3)^2}{s_2}
\]
has a unique minimum in the point $\bar x$ that lies on the geodesic through $x_1$ and $x_3$ and satisfies $d(x_1,\bar x)/s_1=d(x_3,\bar x)/s_2$. We choose exponential coordinates $\xi=\log_{\bar x}x_2$ and expand $f$ in the point $\bar x$. This leads to the following expressions for $f$:
\begin{align*}
  f(x_2)&=\frac{d(x_1,x_3)^2}{s_1+s_2} + \left(\frac1{s_1}+\frac1{s_2}\right)G_{ij}(s_1,s_2,x_1,x_2(\xi),x_3)\xi^i\xi^j\\
  &=\frac{d(x_1,x_3)^2}{s_1+s_2} + \left(\frac1{s_1}+\frac1{s_2}\right) G_{ij}(s_1,s_2,x_1,\bar x,x_3)\xi^i\xi^j\\
  &\quad + G_{ijk}(s_1,s_2,x_1,x_2(\xi),x_3)\xi^i\xi^j\xi^k
\end{align*}
for smooth functions $G_{ij}$ and $G_{ijk}$. The matrix $G_{ij}(s_1,s_2,x_1,x_2,x_3)$ defined and bounded on $([0,1]^2\smallsetminus\{0\})\times K^3$ is positive definite and there is a constant $C>0$ so that the smallest eigenvalue of the matrix is greater than $C$ for all $x_1,x_2,x_3\in K$ and $s_1,s_2\in[0,1]$. Furthermore $G_{ij}$ is homogeneous of degree $0$ in $s_1,s_2$ and we have $G_{ij}(s_1,s_2,x_1,\bar x,x_3)\to\delta_{ij}$ for $|x_1-x_3|\to 0$.
\end{lemma}

\begin{proof}
If $x_2$ is not on the geodesic between $x_1$ and $x_3$, it is easy to see that there is always a point on the geodesic for which one term of $f$ has the same value and the other one is smaller. For $x_2$ on the geodesic we have $d(x_1,x_2)+d(x_2,x_3)=d(x_1,x_3)$ from which $d(x_1,\bar x)/s_1=d(x_3,\bar x)/s_2$ follows. The critical point $\bar x$ of the smooth function $s_1s_2f(x_2)$ is a smooth function of $x_1$, $x_3$, $s_1$, $s_2$, homogeneous of degree zero in $s_1,s_2$, as long as the Hessian is nondegenerate, which is the case if $K$ is small enough and $s_1/(s_1+s_2)\in (-\varepsilon_1,1+\varepsilon_1)$ for some $\varepsilon_1>0$. The expansion of $f$ is just the Taylor expansion (with remainder) in the point $x_2=\bar x$. This gives
\begin{eqnarray*}\lefteqn{
G_{ij}(s_1,s_2,x_1,x_2(\xi),x_3)}\\
&&=\frac1{s_1+s_2}\int_0^1(1-u)\frac{\partial^2}{\partial \eta^i\partial \eta^j}
(s_2d(x_1,\exp_{\bar x}(\eta))^2+s_1d(\exp_{\bar x}(\eta),x_3)^2)|_{\eta=u\xi}du.
\end{eqnarray*}
{}From this expression we see that $G_{ij}$ is homogeneous in $s$ and is smooth for $s_1/(s_1+s_2)\in(-\varepsilon_1,1+\varepsilon_1)$, $x-1,x_3\in K$. In particular it is a bounded continuous function on $([0,1]^2\smallsetminus\{0\})\times K^2$.

For a Euclidean metric it is an application of the law of cosines to show that $G_{ij}=\delta_{ij}$. By rescaling $x_i\mapsto\lambda x_i$ and taking into account that $d(\lambda x_i,\lambda x_j)/\lambda\to|x_i-x_j|$ for $\lambda\to0$, we see that $G_{ij}(s_1,s_2,x_1,x_2,x_3)\to\delta_{ij}$ if $|x_1-x_2|+|x_2-x_3|\to0$. Therefore also $G_{ij}(s_1,s_2,x_1,\bar x,x_3)\to\delta_{ij}$ for $|x_1-x_3|\to0$. As $K$ is small, we are still close to the Euclidean case and therefore $G_{ij}-\delta_{ij}$ is small, from which the existence of $C$ follows.
\end{proof}

\begin{lemma}\label{lem:jlo:comp3}(Asymptotic expansion under the integral) We write $[f(t)]_t$ for the asymptotic expansion of the function $f$ in the variable $t$ in $t=0$. In the following cases we are allowed to interchange the asymptotic expansion and the integration:
\begin{itemize}
\item[(i)] Let $f:[0,1]^{p+1}\setminus\{0\}\to\C$ be a smooth function and assume that there is an $n\in\N$ so that $F(s,t):=t^nf(st)$ can be continued to a function in $C^\infty(\Delta_p\times[0,1])$\footnote{By "$C^\infty$ on a closed set" we mean that every derivative exists in the interior and extends continuously to the boundary.}. Then
\[
  \Big[\int_{\Delta_p}f(st)ds\Big]_t = \int_{\Delta_p}[f(st)]_t\,ds\,.
\]
\item[(ii)] Let $K,G_{ij},G_{ijk},\bar x, x_2(\xi)$ be as in Lemma \ref{lem:jlo:comp2}. Let $H:\R^{2n}\to\C$ be a smooth function with support in a small neighbourhood of the origin. We abbreviate $G_{ij}:=G_{ij}(s_1,s_2,x_1,\bar x,x_3)$, $G_{ij}(\xi):=G_{ij}(s_1,s_2,x_1,x_2(\xi),x_3)$ and $G_{ijk}(\xi):=G_{ijk}(s_1,s_2,x_1,x_2(\xi),x_3)$. Then
\begin{align*}
  &\Big[\int_{\R^{2n}} H(\sqrt{t}\xi)e^{-aG_{ij} (\sqrt{t}\xi)\xi^i\xi^j}\,d\xi\Big]_{\sqrt{t}}\\
  &= \int_{\R^{2n}} [H(\sqrt{t}\xi)e^{-aG_{ijk}(\sqrt{t}\xi)\xi^i\xi^j\xi^k\sqrt{t}}]_{\sqrt{t}} e^{-aG_{ij}\xi^i\xi^j}\,d\xi\,.
\end{align*}
where $a$ is any positive constant.

\end{itemize}
\end{lemma}
\begin{proof}
(i) As $[f(st)]_t=t^{-n}[F(st)]_t$, it suffices to show that we can interchange the integral and the asymptotic expansion for $F$. Because $F$ is smooth, its asymptotic expansion is given by the Taylor series and we have to show that in the following expression the limit and the integral are interchangeable:
\[
  \lim_{t\to0}\int_{\Delta_p}\frac{F(s,t)-\sum_{k=0}^\ell\partial_t^k F(s,t)t^k/k!}{t^{\ell+1}}ds\,.
\]
This is true because the integrand is dominated by $\sup\limits_{t\in[0,1]}|\partial_t^{\ell+1}F(s,t)|/(\ell+1)!$.

(ii) As in part (i), we consider the remainder of the Taylor expansion:
\begin{align*}
  & \frac{(\partial/\partial\sqrt{t})^m}{m!} H(\sqrt{t}\xi)e^{-aG_{ij}(\sqrt{t}\xi)\xi^i\xi^j}\\
  & \quad = \sum_{|\alpha|=m}\xi^\alpha\frac{\partial_\eta^\alpha}{\alpha!} H(\eta)e^{-aG_{ijk}(\eta)\eta^i\xi^j\xi^k}\Big|_{\eta=\sqrt{t}\xi} e^{-aG_{ij}\xi^i\xi^j}\,.
\end{align*}
As $H$ has compact support, we can estimate this by
\[
  \|H(\eta)G_{ijk}(\eta)\eta^i\|_m(1+\|\xi\|^2)^m \xi^\alpha e^{-aG_{ij}\xi^i\xi^j-aG_{ijk}\xi^i\xi^j\xi^k\sqrt{t}}\,.
\]
According to Lemma \ref{lem:jlo:comp2}, there is a constant $C$, so that
\[
  G_{ij}\xi^i\xi^j+G_{ijk}(\sqrt{t}\xi)\xi^i\xi^j\xi^k\sqrt{t} = G_{ij}(\sqrt{t}\xi)\xi^i\xi^j>C\|\xi\|^2\,,
\]
for all $\xi$ such that $\sqrt t\xi$ is in the support of $H$. Thus it follows again by the dominated convergence theorem that the asymptotic expansion and the integral commute.
\end{proof}

\begin{proof}[Proof of Proposition \ref{prop:jlo:comp_power}.] We write Latin letters for indices in $\N_0$ and Greek letters for multiindices in $\N_0^{2n}$.

We consider again the function $f(s):=\Str(D_0 K_{s_0}^N [\barpart^*,D_1]K_{s_1}^N \cdots [\barpart^*,D_p] K_{s_p}^N)$. As we are going to show, the asymptotic expansion of $f(st)$ w.r.t. $t$ in $t=0$ exists, has lowest order $-n-\lfloor\frac d2\rfloor$ and the coefficients are smooth functions of $s\in\Delta_p$. Therefore the function $F(s,t):=t^{n+\lfloor\frac d2\rfloor}f(st)$ is smooth\footnote{From Proposition \ref{prop:jlo:fsmooth} follows that $F$ is smooth for $(s,t)\in\Delta_p\times(0,1]$. The existence of the asymptotic expansion shows that its derivatives can be continued to $t=0$. Hence $F\in C^\infty(\Delta_p\times[0,1])$} and we can apply Lemma \ref{lem:jlo:comp3} (i):
\[
  \Psi_p(D_0,\dots,D_p) = (-1)^{\frac{p(p+1)}2} \int_{\Delta_p} [t^pf(st)]_- ds\,.
\]
To compute the asymptotic expansion of $f(st)$, we consider the kernel
\[
  (D_0 K_{s_0}^N [\barpart^*,D_1]K_{s_1}^N \cdots [\barpart^*,D_p] K_{s_p}^N)(x_0,x_{p+1})\,.
\]
Recall that $D_0$ has compact support $K\subset U\subset X$ where $U$ is an open set (see also Proposition \ref{p-Lili}). As $K_{s_i}^N(x_i,x_{i+1})$ vanishes for $d(x_i,x_{i+1})>\varepsilon$, there is a $\varepsilon>0$ so that $(p+1)\varepsilon$ is smaller than the geodesic distance between $K$ and $X\setminus U$. Then in the above kernel only the values of terms inside a compact subset $K_\varepsilon$ of $U$ play a role and therefore it is well defined. We assume that $K_\varepsilon$ is small enough to apply Lemma \ref{lem:jlo:comp2}.

We want to ``move'' the operators $[\barpart^*,D_i]$ to the left. First just consider a term $K_{s_1}^N D K_{s_2}^N$. We may assume that $D$ in local coordinates has the form $\rho(x)\partial^\alpha$ where $\supp\rho\subset K_\varepsilon$. Explicitly, the above term is given by the integral
\[
  \int\limits_X \sum_{0\le i,j\le N} s_1^is_2^j\Psi_i(x_1,x_2)\frac{e^{-d(x_1,x_2)^2/s_1}}{(\pi s_1)^n} \rho(x_2)\partial^\alpha_{x_2}\left(\Psi_j(x_2,x_3) \frac{e^{-d(x_2,x_3)^2/s_2}}{(\pi s_2)^n}\right) |dx_2|_g\,.
\]
We write $|dx_2|_g = \sigma(x_2)dx_2$ and integrate by parts to bring the $\partial^\alpha_{x_2}$-operator to the left. Then we make repeatedly use of Lemma \ref{lem:jlo:comp1} to "replace" the $x_2$-derivatives by $x_1$-derivatives, i.e. we use an identity of the form
\[
  \partial^{\alpha}_{x_2}e^{-d(x_1,x_2)^2/s_1} = \sum_{\beta+\gamma=\alpha} h_{\beta,\gamma}(x_1,x_2) \partial^{\gamma}_{x_1} e^{-d(x_1,x_2)^2/s_1}\,,
\]
which holds for some smooth functions $h_{\beta,\gamma}$. Writing down again the integral, we find an expression of the form
\[
  \int_X \sum_{0\le i,j\le N} \sum_{\alpha'\le\alpha} s_1^is_2^j H_{i,j,\alpha'}(x_1,x_2,x_3)\partial^{\alpha'}_{x_1} \frac{e^{-d(x_1,x_2)^2/s_1-d(x_2,x_3)^2/s_2}}{(\pi s_1)^n(\pi s_2)^n} dx_2
\]
where $H_{i,j,\alpha'}$ are smooth functions.
If we apply the above procedure to shift all derivatives in the expression $D_0K^N_{s_0}\dots [\barpart^*,D_p]K^N_{s_p}$ to the left, we get
\begin{equation}\label{F:jlo:as2}
  \int_{X^p} \sum_{|\gamma|\le N} \sum_{|\alpha|\le d} s^\gamma H_{\gamma,\alpha}(x_0,\dots,x_p)\partial^{\alpha}_{x_0} \frac{e^{-\sum_{j=0}^p d(x_j,x_{j+1})^2/s_j}} {(\pi s_0)^n\dots(\pi s_p)^n} dx_1\dots dx_p\,.
\end{equation}
We omitted the terms for which $|\gamma|:=\sum_{j=1}^{p-1}\gamma_j>N$, but we will see later that they would only produce (irrelevant) terms of higher order in $t$. We rewrite the exponent in the above expression using Lemma \ref{lem:jlo:comp2} repeatedly:
\begin{eqnarray*}\lefteqn{
  \sum_{j=0}^{p}\frac{d(x_j,x_{j+1})^2}{s_j}= \frac{d(x_0,x_{p+1})^2}{s_0+\dots+s_p}
}\\
 &&+ \sum_{\ell=1}^p \left(\frac1{s_0+\dots+s_{\ell-1}}+\frac1{s_\ell}\right) G_{ij}(s_0+\dots+s_{\ell-1},s_\ell,x_0,x_\ell,x_{\ell+1})\xi^i_\ell\xi^j_\ell,
\end{eqnarray*}
where $\xi_\ell=\ln_{\bar x_\ell}x_\ell$, $\bar x_\ell= \bar x_\ell(x_{\ell-1},x_{\ell+1})$. Now we change to the variables $\xi^i$ in the integral and rescale $\xi^i\mapsto\sqrt{t}\xi^i$ as well as $s_i\mapsto ts_i$ so that $(s_0,\dots,s_p)\in\Delta_p$. We temporarily forget the last term in the exponent and suppress the arguments of $G_{ij}$:
\begin{equation}\label{F:jlo:as3}
  t^{p-n}\!\!\!\!\!\int\limits_{(T_{\bar x}X)^p} \sum_{\abtop{|\gamma|\le N}{|\alpha|\le d}}\! s^\gamma H_{\gamma,\alpha}(x_0,\dots,x_{p+1})\partial^{\alpha}_{x_0} \frac{e^{-\sum_{\ell=1}^p\left(\frac1{s_0+\dots+s_{\ell-1}}+\frac1{s_\ell}\right) G_{ij}\xi^i_\ell\xi^j_\ell}}{(\pi s_0)^n\dots(\pi s_p)^n} d\xi_1\dots d\xi_p
\end{equation}
where the Jacobi determinant has been absorbed in $H_{\gamma,\alpha}$. Due to Lemma \ref{lem:jlo:comp3} ii) we are allowed to expand asymptotically w.r.t. $\sqrt{t}$ under the integral. Keep in mind that $x_\ell=\exp_{\bar x_\ell}(\sqrt{t}\xi_\ell)$ so that the arguments of $H_{\gamma,\alpha}$ as well as of $G_{ij}$ depend on $\sqrt{t}$. In the expansion of the exponent, there will be singular terms in $s$, namely powers of the factor
\[
 \frac1{s_0+\dots+s_{\ell-1}}+\frac1{s_\ell}=\frac{s_0+\dots+s_\ell}{(s_0+\dots+s_{\ell-1})s_\ell}\,,
\]
but as these factors only appear paired with $\xi^i_\ell\xi^j_\ell$, the singularities cancel as we see in the following computation. After the expansion we have to compute integrals of the form
\[
  \int_{T_{\bar x}X} \xi_\ell^\beta e^{-\frac{s_0+\dots+s_\ell}{(s_0+\dots+s_{\ell-1})s_\ell} G_{ij}\xi^i_\ell\xi^j_\ell }d\xi_\ell = C_\beta(s,x)\left(\frac{(s_0+\dots+s_{\ell-1})s_\ell}{(s_0+\dots+s_\ell)}\right)^{\frac{|\beta|}2+n}
\]
where $C_\beta(s,x)$ is a smooth function, homogeneous of degree 0 in $s$, vanishing unless $|\beta|=\sum\beta_i$ is even.
Terms with $|\beta|$ even correspond to even terms in the asymptotic expansion in powers of $\sqrt t$. Therefore we actually have an asymptotic series in $t$. 

We repeat the above steps for $\xi_2,\dots,\xi_p$. As
\[
  \prod_{\ell=1}^p\frac{(s_0+s_1+\dots+s_{\ell-1})s_\ell}{s_0+s_1+\dots+s_\ell} = \frac{s_0s_1\dots s_p}{s_0+s_1+\dots+s_p}\,,
\]
the singularities from the denominator in equation (\ref{F:jlo:as3}) disappear, and we get
\begin{align*}
  &(D_0K^N_{ts_0}\cdots D_pK^N_{ts_p})(x_0,x_{p+1})\\
  &=t^{p-n}\sum_{|\alpha|\le d}\sum_{k=0}^N t^kf_k(s,x_0,x_{p+1})\partial_{x_0}^\alpha e^{-\frac{d(x_0,x_{p+1})^2}{t}}+\cO(t^{p-n+N+1})
\end{align*}
for smooth functions $f_k:\Delta_p\times K\times K_\varepsilon\to\C$. Remember that
\[
  f(s,t)=\int_K (D_0K^N_{ts_0}\cdots D_pK^N_{ts_p})(x_0,x_0) dx_0\,.
\]
The integral over $x_0\in K$ and the asymptotic expansion commute for the same reason as in Lemma \ref{lem:jlo:comp3}. We see in the above formula that the negative powers in $t$ are only produced by the derivative $\partial_{x_0}^\alpha$. As $\lim\limits_{x_p\to x_0} \partial_{x_0}^\alpha d(x_0,x_p)^2 =0$ for $|\alpha|=1$, we need at least two derivatives to get a negative power in $t$. Thus the negative power is at most $\lfloor\frac{|\alpha|}{2}\rfloor$.

In formula (\ref{F:jlo:as2}), the coefficient functions of the operators $D_0,D_1,\dots,D_p$ have been absorbed in the function $H_{\gamma,\alpha}$. It is easy to check that they enter linearly and with derivatives of order at most $d$, which is the sum of the degrees of the differential operators, in this function. After formula \eqref{F:jlo:as3} when we do the expansion, we get an additional derivative for every order of $\sqrt{t}$. Therefore the coefficients of $\Psi_p$ only depend on finitely many derivatives of the operators $D_0,D_1,\dots,D_p$ restricted to the compact set $K_\varepsilon$ that was mentioned in the beginning of the proof. This means that we can estimate $\Psi_p$ by a product of $C^k$-norms over the compact set $K_\varepsilon$ of the operators $D_i$. As the operators $D_1,\dots,D_p$ are holomorphic and actually defined on an open set containing $K_\varepsilon$, we can use the Cauchy integral formula to estimate their $C^k$-norms by the sup-norms over a compact set that is slightly bigger than $K_\varepsilon$. This shows that $\Psi_p$ is continuous in the operators $D_0,\dots,D_p$.
\end{proof}

\begin{remark}\label{rem:jlo:replacebyfullhk}
If we replace in the formula for $\Psi_p$ one of the approximated heat kernels $k^N$ by the exact heat kernel $k$, the $dx$-integral still is over a compact set. Thus we can choose $\varepsilon$ small enough so that the formula for $\Psi_p$ is still well defined. From the above proof it is also clear that this procedure doesn't change the value of $\Psi_p$.
\end{remark}


\begin{bibdiv}
\begin{biblist}
\bib{BS}{article}{
   author={Be{\u\i}linson, A. A.},
   author={Schechtman, V. V.},
   title={Determinant bundles and Virasoro algebras},
   journal={Comm. Math. Phys.},
   volume={118},
   date={1988},
   number={4},
   pages={651--701},
   issn={0010-3616},
   review={\MR{962493 (90m:32048)}},
}
\bib{BGV}{book}{
   author={Berline, Nicole},
   author={Getzler, Ezra},
   author={Vergne, Mich{\`e}le},
   title={Heat kernels and Dirac operators},
   series={Grundlehren der Mathematischen Wissen\-schaften [Fundamental
   Principles of Mathematical Sciences]},
   volume={298},
   publisher={Springer-Verlag},
   place={Berlin},
   date={1992},
   pages={viii+369},
   isbn={3-540-53340-0},
   review={\MR{1215720 (94e:58130)}},
}

\bib{BR}{article}{
   author={Bern{\v{s}}te{\u\i}n, I. N.},
   author={Rosenfel{\cprime}d, B. I.},
   title={Homogeneous spaces of infinite-dimensional Lie algebras and the
   characteristic classes of foliations},
   language={Russian},
   journal={Uspehi Mat. Nauk},
   volume={28},
   date={1973},
   number={4(172)},
   pages={103--138},
   issn={0042-1316},
   review={\MR{0415633 (54 \#3714)}},
   translation={
      journal={Russian Math. Surveys},
      volume={28},
      date={1973},
      pages={107--142},
      number={4},
   },
}

\bib{BG}{article}{
   author={Brylinski, Jean-Luc},
   author={Getzler, Ezra},
   title={The homology of algebras of pseudodifferential symbols and the
   noncommutative residue},
   journal={$K$-Theory},
   volume={1},
   date={1987},
   number={4},
   pages={385--403},
   issn={0920-3036},
   review={\MR{920951 (89j:58135)}},
}
\bib{Connes}{article}{
   author={Connes, Alain},
   title={Noncommutative differential geometry},
   journal={Inst. Hautes \'Etudes Sci. Publ. Math.},
   number={62},
   date={1985},
   pages={257--360},
   issn={0073-8301},
   review={\MR{823176 (87i:58162)}},
}

\bib{FLS}{article}{
    title = {{Riemann-Roch-Hirzebruch theorem and Topological Quantum
        Mechanics}},
    author = {Feigin, Boris},
    author = {Losev, Andrey},
    author = {Shoikhet, Boris},
    eprint = {arXiv:math.QA/0401400},
}
\bib{FT}{article}{
   author={Fe{\u\i}gin, B. L.},
   author={Tsygan, B. L.},
   title={Riemann-Roch theorem and Lie algebra cohomology. I},
   booktitle={Proceedings of the Winter School on Geometry and Physics
   (Srn\'\i, 1988)},
   journal={Rend. Circ. Mat. Palermo (2) Suppl.},
   number={21},
   date={1989},
   pages={15--52},
   review={\MR{1009564 (90k:17041)}},
}
\bib{FFS}{article}{
   author={Feigin, Boris},
   author={Felder, Giovanni},
   author={Shoikhet, Boris},
   title={Hochschild cohomology of the Weyl algebra and traces in
   deformation quantization},
   journal={Duke Math. J.},
   volume={127},
   date={2005},
   number={3},
   pages={487--517},
   issn={0012-7094},
   review={\MR{2132867 (2006j:16015)}},
}
\bib{G}{article}{
   author={Gel{\cprime}fand, I. M.},
   title={The cohomology of infinite dimensional Lie algebras: some
   questions of integral geometry},
   conference={
      title={Actes du Congr\`es International des Math\'ematiciens},
      address={Nice},
      date={1970},
   },
   book={
      publisher={Gauthier-Villars},
      place={Paris},
   },
   date={1971},
   pages={95--111},
   review={\MR{0440631 (55 \#13505)}},
}

\bib{GK}{article}{
   author={Gel{\cprime}fand, I. M.},
   author={Ka{\v{z}}dan, D. A.},
   title={Certain questions of differential geometry and the computation of
   the cohomologies of the Lie algebras of vector fields},
   language={Russian},
   journal={Dokl. Akad. Nauk SSSR},
   volume={200},
   date={1971},
   pages={269--272},
   issn={0002-3264},
   review={\MR{0287566 (44 \#4770)}},
}

\bib{GKF}{article}{
   author={Gel{\cprime}fand, I. M.},
   author={Ka{\v{z}}dan, D. A.},
   author={Fuks, D. B.},
   title={Actions of infinite-dimensional Lie algebras},
   language={Russian},
   journal={Funkcional. Anal. i Prilo\v zen.},
   volume={6},
   date={1972},
   number={1},
   pages={10--15},
   issn={0374-1990},
   review={\MR{0301767 (46 \#922)}},
}


\bib{JLO}{article}{
   author={Jaffe, Arthur},
   author={Lesniewski, Andrzej},
   author={Osterwalder, Konrad},
   title={Quantum $K$-theory. I. The Chern character},
   journal={Comm. Math. Phys.},
   volume={118},
   date={1988},
   number={1},
   pages={1--14},
   issn={0010-3616},
   review={\MR{954672 (90a:58170)}},
}

\bib{Lefschetz}{book}{
   author={Lefschetz, Solomon},
   title={Introduction to Topology},
   series={Princeton Mathematical Series, vol. 11},
   publisher={Princeton University Press},
   place={Princeton, N. J.},
   date={1949},
   pages={viii+218},
   review={\MR{0031708 (11,193e)}},
}
\bib{Lysov}{article}{
   author={Lysov, Vyacheslav},
   title={Anticommutativity equations in topological quantum
   mechanics},
   journal={Pisma Zh.Eksp.Teor.Fiz.},
   volume={76},
   date={2002},
   pages={855--858},
   translation={
   journal={JETP Lett.},
   volume={76},
   date={2002},
   pages={724--727}},
}
\bib{Loday}{book}{
   author={Loday, Jean-Louis},
   title={Cyclic homology},
   series={Grundlehren der Mathematischen Wissenschaften [Fundamental
   Principles of Mathematical Sciences]},
   volume={301},
   edition={2},
   publisher={Springer-Verlag},
   place={Berlin},
   date={1998},
   pages={xx+513},
   isbn={3-540-63074-0},
   review={\MR{1600246 (98h:16014)}},
}
\bib{Maclane}{book}{
   author={MacLane, Saunders},
   title={Homology},
   edition={1},
   note={Die Grundlehren der mathematischen Wissenschaften, Band 114},
   publisher={Springer-Verlag},
   place={Berlin},
   date={1967},
   pages={x+422},
   review={\MR{0349792 (50 \#2285)}},
}
\bib{NT1}{article}{
   author={Nest, Ryszard},
   author={Tsygan, Boris},
   title={Algebraic index theorem},
   journal={Comm. Math. Phys.},
   volume={172},
   date={1995},
   number={2},
   pages={223--262},
   issn={0010-3616},
   review={\MR{1350407 (96j:58163b)}},
}

\bib{R}{article}{
    title = {Some notes on the Feigin--Losev--Shoikhet integral conjecture},
    author = {Ramadoss, Ajay},
    eprint = {arXiv:math.QA/0612298},
}

\bib{S}{article}{
   author={Schechtman, V. V.},
   title={Riemann-Roch theorem after D. Toledo and Y.-L. Tong},
   booktitle={Proceedings of the Winter School on Geometry and Physics
   (Srn\'\i, 1988)},
   journal={Rend. Circ. Mat. Palermo (2) Suppl.},
   number={21},
   date={1989},
   pages={53--81},
   review={\MR{1009565 (91b:58242)}},
}

\bib{Treves}{book}{
   author={Tr{\`e}ves, Fran{\c{c}}ois},
   title={Topological vector spaces, distributions and kernels},
   publisher={Academic Press},
   place={New York},
   date={1967},
   pages={xvi+624},
   review={\MR{0225131 (37 \#726)}},
}

\bib{Wodzicki}{article}{
   author={Wodzicki, Mariusz},
   title={Cyclic homology of differential operators},
   journal={Duke Math. J.},
   volume={54},
   date={1987},
   number={2},
   pages={641--647},
   issn={0012-7094},
   review={\MR{899408 (88k:32035)}},
}

\end{biblist}
\end{bibdiv}
\end{document}